\magnification\magstep1

\def\d{{\rm deg}}

\def\cA{{\cal A}}
\def\cB{{\cal B}}
\def\cC{{\cal C}}
\def\cD{{\cal D}}
\def\cE{{\cal E}}
\def\cF{{\cal F}}
\def\cH{{\cal H}}
\def\cI{{\cal I}}
\def\cL{{\cal L}}

\def\cO{{\cal O}}
\def\cP{{\cal P}}

\def\cS{{\cal S}}
\def\cT{{\cal T}}

\def\C{{\bf C}}

\def\P{{\bf P}}

\centerline{\bf On the Degree Growth of  Birational Mappings}

\centerline{\bf in Higher Dimension}
\medskip
\centerline{Eric Bedford and Kyounghee Kim}
\bigskip
\centerline{\bf \S0.  Introduction}
\medskip
\noindent
Let $f:{\bf C}^d\to{\bf C}^d$ be a birational map.  The problem of
determining the behavior of the iterates $f^n=f\circ\cdots\circ f$ is very
interesting but not well understood.  A basic property of a
rational map is its degree (see \S1 for the definition).  Another
quantity is the dynamical degree 
$$\delta(f)=\lim_{n\to\infty}({\rm deg}(f^n))^{1\over n}$$
which is invariant under birational self-maps of $\C^d$ (see [BV]).   It has been called ``complexity'' by some physicists (see [BM] and
[AABM]), and its logarithm has been called ``algebraic entropy'' in [BV].

One aspect of a birational map is that its birational conjugacy class
does not have a well-defined ``domain.'' Namely, if $f:X\to X$ is a (birational)
dynamical system, then any birational equivalence $h:X\to\tilde X$ will convert
$f$ to a (birational) dynamical system on $\tilde f = h\circ f\circ
h^{-1}:\tilde X\to\tilde X$.  This may serve as a significant change of the
presentation of $f$, since the cohomology groups of $X$ and $\tilde X$ may have
different dimensions.   There is a well-defined pull-back map on cohomology
$f^{*}:H^{1,1}(X)\to H^{1,1}(X)$.   Thus we can define $\delta(f)$ more generally as $\lim_{n\to\infty}||(f^n)^*||^{1/n}$, which is the exponential rate of growth of the action of $f^n$ on $H^{1,1}$.  Dinh and Sibony [DS] showed that this more general $\delta(f)$ is birationally invariant.

The passage to cohomology may or may not be compatible
with the dynamical system, depending on whether $(f^{*})^{n}=(f^{n})^{*}$
holds on
$H^{1,1}$.  In case this holds, we say that
$f$ is {\it $H^{1,1}$-regular}, or simply {\it 1-regular}.   And if $f$ is 1-regular, then it follows that
$\delta$ is the spectral radius of $f^{*}$.  We will pursue the study of
$\delta$ using what might be called a method of regularization.  The first step
of this method is to replace the pair $(f,X)$ by a regular pair $(\tilde
f,\tilde X)$.  This is done by finding a new complex manifold $\tilde X$ which
is birationally equivalent to $X$.  The second step, then, is to determine
$\tilde f^{*}$ and its spectral radius.

In this paper, we show how the method of regularization may be carried out on
certain sub-classes of the family of maps which have the form $f=L\circ J$,
where $L$ is an invertible linear map, and
$J(x_{1},\dots,x_{d})=(x_{1}^{-1},\dots,x_{d}^{-1})$.  Generally speaking, a birational map has subvarieties that are mapped to lower
dimensional sets, and it has lower dimensional subvarieties that are blown up
to sets of higher dimension.  In the case of
$f=L\circ J$, the coordinate hypersurfaces $\{x_{j}=0\}$ are blown down to
points, and the points $e_\ell=[0:\dots:1:\dots:0]$ are blown back up to hypersurfaces.  This
interplay between the blowing down and the blowing up serves as the key
for degree growth.  The essential objects are those orbits of the
following form:
$\{x_{j}=0\}\to *\to\cdots\to*\to e_\ell$.  That is, they start at a hypersurface
which is blown down to a point, and the orbit of this point lands on a point of
indeterminacy which is blown back up to a hypersurface.  We say that such orbits are singular.  In \S4 we define
the class of elementary maps, which are characterized by the property that they are locally biholomorphic at the intermediate points of singular orbits.  The 1-regularization of an elementary map is obtained by blowing up the points of the singular orbits.  The singular orbits of an elementary map can be organized into an orbit list structure
$\cL^{c}$, $\cL^{o}$, which consists of two sets of lists of positive integers.  It is
then shown how $\delta=\delta(\cL^{c},\cL^{o})$ is determined by this list
structure: an expression for $\tilde f^{*}$ is given in \S4, and the
characteristic polynomial is given in the Appendix.  In \S7 we illustrate this method by carrying out the procedure of finding orbit lists for some maps that have appeared in the mathematical physics literature.
 
In general we have $\delta\le\deg(f)$, and the existence of a singular orbit
causes the inequality to be strict (Theorem 4.2).  It is often 
desirable to have some way of estimating $\delta$ without actually computing
it.  To this end, we give a way of comparing
$\delta$ for two maps $f$
and $\hat f$ with list structures $\cL^{c}$, $\cL^{o}$ and
$\hat\cL^{c}$,
$\hat\cL^{o}$, respectively.  In \S5 we give a number of comparison results;
here are two examples.  Theorem 5.1 shows that if $f$ and
$\hat f$ have the same list structures, except that the orbits of $\hat f$ are
longer, then
$\delta(f)\le\delta(\hat f)$.  In Theorem 5.3, we show that adding a complete
orbit list decreases $\delta$.  If we simply add a new singular orbit of length
$M$ to one of the orbit lists, then whether $\delta$ is increased or decreased
depends on the size of $M$.

In \S6 we introduce linear maps $L_p$, which are determined by a permutation $p$.  The rest of  \S6 will be devoted to a consideration of the case where $p=I$ is the identity permmutation.  The maps $f=L_I\circ J$ are the Noetherian maps which were defined in the work [BHM].  Our primary motivation in this section is to show how the work of \S4 applies to give $\deg(f^n)$, $n\ge0$ for these maps.  These same numbers were conjectured in [BHM].

For more general permutations $p$, the map $f_p=L_p\circ J$ can lead to some complicated examples of orbit collision.  The expression {\it orbit collision} refers to the fact that a singular orbit $\{x_j=0\}\to*\to\dots\to\sigma_k\to\dots\to e_\ell$ can contain a point $\sigma_k\in\{x_k=0\}-\cI$ where $f$ is smooth but not locally invertible.  Thus the singular orbit starting  at $\{x_j=0\}$ contains the singular orbit starting at $\{x_k=0\}$ (and possibly others).   In \S8 we define singular chains.  The singular chain structure is used both to construct the 1-regularization $(f_X,X)$  of $(f_p,\P)$ and to write down $f^*_X$.  As before, $\delta(f_p)$ is given by the spectral radius of $f^*_X$.

In this paper we follow up on ideas of Diller and Favre [DF], Guedj [G] and Boukraa, Hassani and Maillard [BHM].  The paper [DF] shows that 1-regularization is possible for all birational maps in dimension two, and this forms the basis for their penetrating analysis.  The possibility of extending this approach to the case of higher
dimension was discussed in [BHM].  Birational maps are more complicated
in higher dimension, however, and [BHM] proposed a family of birational
maps as a model family of mappings to analyze.  In their study of these
mappings, they identify the integrable cases and give a numerical
description of $\delta(f)$ in many cases.  This family was chosen in part
because it resembles maps that arise in mathematical physics (see [AABHM],
[BTR], and [RGMR]).  Our analysis of elementary maps evolved from an effort to understand questions posed in [BHM].

The contents of this paper are as follows.  In \S1 we assemble some basic
concepts concerning rational maps, and we formulate the property (1.1) which we
use in constructing 1-regularizations.  In \S2 we recall the basic relationship
between $(f^*)^n$ on cohomology and the degree of $f^n$.  In \S3 we show how to
1-regularize the basic mapping $J$.  Then in \S4 we apply this to elementary
maps and show how to compute the induced mapping $f^*$ on $H^{1,1}$.  At this
stage, we may compute the characteristic polynomial $\chi_f$ of $f^*$.  The
actual computation is deferred to Appendix A so that our discussion of
$\delta$ is not interrupted. 
\S5 shows how to use to formula for $\chi_f$ to obtain comparison theorems for
$\delta$.  In \S6 we introduce the family of linear transformations $L_p$,
which depend on a permutation $p$, and we give the degree growth for Noetherian maps.  In \S7 we analyze some mappings that have
appeared in the mathematical physics literature.  In \S8 we give the general method for 1-regularization of permutation mappings.

\bigskip
\centerline{\bf \S1.  Polynomial and Rational Maps}
\medskip\noindent
We will review some of the basic properties of the dynamical degree of
birational maps.  We refer the reader to [RS] and [S] for further details
for birational maps of ${\bf P}^d$ and to [DS] for the case of more general
manifolds.

 A polynomial is a finite sum
$p(x)=\sum c_{i_1\dots i_d}x_1^{i_1}\cdots x_d^{i_d}$.  We define the
degree of a monomial as
$\d(x_1^{i_1}\cdots x_d^{i_d})=i_1+\cdots+i_d$, and we define $\d(p)$ to be
the maximum of
$i_1+\cdots+ i_d$ for all nonzero coefficients
$c_{i_1\dots i_d}$.  A mapping $p=(p_1,\dots,p_d)$  is
said to be a polynomial mapping if each coordinate function is polynomial. 
Similarly, a mapping $f=(f_1,\dots,f_d)$ is said to be rational if each
coordinate function is rational, i.e., each $f_j=p_j/q_j$ is a quotient
of polynomials.  

We will use
complex projective space $\P^d$ as a compactification of $\C^d$.  Recall
that 
$$\P^d=\{[x_0:\cdots:x_d]:x_j\in\C,{\rm\ and\ the\ }x_j{\rm\ are\ not\ all\
}0\},$$
where the notation $[x_0:\cdots:x_d]$ denotes homogeneous coordinates,
i.e., $[x_0:\cdots:x_d]=[\zeta x_0:\cdots:\zeta x_d]$ for all nonzero
$\zeta\in\C$. 

A rational map $f=(p_1/q_1,\dots,p_d/q_d)$ of $\C^d$ induces a (partially
defined) holomorphic map 
$\hat f=[\hat f_0:\cdots:\hat f_d]$ of projective space.  Let us describe
how to obtain
$\hat f$ from $f$. 
We add the variable $x_0$ and convert $f$ to a homogeneous function
$$\tilde f: [x_0:\cdots:x_d]\mapsto\left[1: {p_1(x_1/x_0,\dots,x_d/x_0)\over
q_1(x_1/x_0,\dots,x_d/x_0)}:\cdots: {p_d(x_1/x_0,\dots,x_d/x_0)\over
q_d(x_1/x_0,\dots,x_d/x_0)}\right]$$
$$=\left[ 1:{\tilde p_1(x_0,\dots,x_d)\over
\tilde q_1(x_0,\dots,x_d)}:\cdots: {\tilde p_d(x_0,\dots,x_d)\over
\tilde q_d(x_0,\dots,x_d)}\right].$$
Note that the first line
is homogeneous of degree zero, and thus $\tilde p_j$ and $\tilde
q_j$ are homogeneous polynomials with $\d(\tilde p_j)=\d(\tilde q_j)$.  We
passed from the first equation to the second by multiplying the numerators
and denominators by powers of $x_0$.  
Let $Q=\tilde q_1\cdots \tilde q_d$.  On the dense set where $Q\ne0$, we do
not change $\tilde f$ if we multiply by $Q$.  Thus $\tilde f$ is also given
by the map $x\mapsto[Q(x):Q(x)\tilde p_1(x)/\tilde q_1(x):\cdots: \tilde
p_d(x)/\tilde q_d(x)]$.   Thus we have represented $f$ by a polynomial
mapping to projective space. To obtain
$\hat f=[\hat f_0:\cdots:\hat f_d]$, we divide out the greatest
(polynomial) factor.  After this is done, there is no polynomial that
divides all the $\hat f_j$, and we define $\d(f):=\d(\hat f_0)=\cdots=\d(\hat f_d)$.
We may also take the $n$-fold composition $f^n=f\circ\cdots\circ f$, and
perform the same passage to a map $\widehat{f^n}$ on projective space. 
Since we may divide out a (possibly larger) common factor at the end, it is
evident that
$$\d(f^n)\le (\d(f))^n.$$

The indeterminacy locus is the set 
$$\cI( f)=\{x\in\P^d:\hat f_0(x)=\cdots=
\hat f_d(x)=0\}.$$  It is evident that $\hat f$ defines a holomorphic mapping
of $\P^d-\cI(f)$ to $\P^d$.  We denote this simply by $f$; we reserve the notation $\hat f$ for the multiple-valued mapping which will be defined below.  In fact, ${\bf
P}^d-\cI$ is the maximal domain on which $f$ can be extended to be
analytic.  For if 
$\hat a\in \cI(\hat f)$, there is no neighborhood $\omega$ of $\hat a$ such that
$f(\omega-\cI)$ is relatively compact in any of the affine coordinate charts
$U_j=\{x_j\ne0\}$.  For in this case, at least one of the coordinate functions, say
$\hat f_1$, has no common factor with $\hat f_j$.  Thus $(\omega-\{\hat
f_j=0\})\ni x\mapsto\hat f_1(x)/\hat f_j(x)$ will take on all complex values.  For
$a\in{\bf P}^d$, let us define the cluster set $Cl_f(a)$ to be the set of all
limits of $f(a')$ for $a'\in {\bf P}^d-\cI$, $a'\to a$.  The cluster set is
connected and compact.  By the arguments above, it follows that $Cl_f(a)$ contains
more than one point exactly when $a\in\cI$.

Let us consider the graph of the restriction of $f$ to ${\bf P}^d-\cI$:
$$\Gamma_f=\{(x,y):x\in{\bf P}^d-\cI, y=f(x)\}.$$
Thus $\Gamma_f$ is a subvariety of $({\bf P}^d-\cI)\times{\bf P}^d$, and by
$\hat\Gamma_f$ we denote the closure of $\Gamma_f$ inside ${\bf P}^d\times{\bf
P}^d$.  Thus $\hat\Gamma_f$ is an algebraic variety.  Let $\pi_1$ (resp. $\pi_2$)
denote the projections of $\hat\Gamma_f$ to the first (resp. second) coordinate. 
Now we may define the multiple-valued mapping $\hat
f(x):=\pi_2\circ\pi_1^{-1}(x)$, and thus $\hat f(x)$ is a subvariety of $X$ for
each $x$.  We have $\hat f(x)=Cl_f(x)$, and the dimension of $\hat f(x)$ is greater
than zero exactly when $x\in\cI$.

A projective manifold $X$ is said to be rational if it is birationally equivalent to
${\bf P}^d$.  The discussion above applies to rational manifolds.  Let $\cT(X)$
denote the set of positive, closed currents on $X$ of bidegree (1,1).  One of the
well-known properties of a positive, closed (1,1)-current $T$ is that it has a local
potential $p$ and can be written locally as $T = dd^cp$.  Following Guedj [G] we
use local potentials to define the induced pull-back map
$\Phi_f:\cT(X)\to\cT(X)$.  Namely, if $T\in\cT(X)$, and if $x_0\in X-\cI$, then $T$ has a local potential $p$ in a neighborhood for $f(x_0)$, and we define the pullback
$f^*T:=dd^c(p\circ f)$ in a neighborhood of $x_0$.  This yields a well-defined,
positive, closed (1,1)-current on the set $X-\cI$.  Now by [HP], the set $\cI$,
being a subvariety of codimension at least 2, is a ``removable singularity'' for a
positive, closed (1,1) current.  This means two things.  First, the current
$f^*T$ has finite total mass, so it may be considered to be a (1,1)-form whose
coefficients are (complex, signed) measures with finite total mass.  This allows us
to define $\Phi_f(T):=\widetilde{f^*T}$ as the current obtained by extending these measures ``by zero'' to $X$, i.e. by assigning zero mass to the set
$\cI$.  Second, the current $\Phi_f(T)$ is closed.  Thus
$\Phi_f(T)\in \cT(X)$.

The currents we will use are currents of integration.  Specifically, if $V$ is a
subvariety of pure codimension 1, then we define the current of integration $[V]$
as an element of the dual space to the space of smooth $(d-1,d-1)$-forms:
$\varphi\mapsto\langle [V],\varphi\rangle:=\int_V\varphi$.  By a classic theorem of
Lelong, $[V]$ is well-defined and is a positive, closed (1,1)-current.  If $V$ is
defined locally as $\{h=0\}$ for some holomorphic function $h$, then ${1\over 2\pi}\log |h|$ is a local potential for $[V]$, which means that locally, $[V]={1\over
2\pi}dd^c\log|h|$.  The pull-back of the current in this case is simply the
preimage: $\Phi_f([V]) = [(f|_{X-\cI})^{-1}V]$.

An irreducible subvariety $V$ will be said to be {\it exceptional} if
$V-\cI(f)\ne\emptyset$, and if ${\rm dim}(f(V-\cI(f)))<{\rm dim}(V)$.  The exceptional locus for $f$, written $\cE(f)$ is the union of all irreducible exceptional varieties.

We will use the following condition:
$$\eqalign{ &{\rm For\ every\  exceptional\
hypersurface\ }V{\rm \ and\ every\ }n>1,\cr
&{  \rm the\ image\ }\hat f^n(V-\cI){\rm \ has\ codimension\
strictly\ greater\ than\ one.} } \eqno(1.1)$$  
Note that if $V$ is exceptional, then $f(V-\cI)=\hat f(V-\cI)$ will be contained in a subvariety of codimension at least 2.  The only possibility that the codimension could jump to 1 for $n=2$ would come from the action of $\hat f$ on $\cI\cap(f(V-\cI))$.  Thus (1.1) depends on the
behavior of $\hat f$ at certain points of $\cI$.  A related condition, called algebraic
stability, was introduced in  [FS] for maps of $\P^{d}$ and is equivalent to
$\deg(f^{n})=(\deg(f))^{n}$ for all $n\ge0$.
\proclaim Proposition 1.1.  If (1.1) holds, then $(\Phi_f)^n=\Phi_{f^n}$.

\noindent{\it Proof.}  
Let us fix $T\in\cT(X)$.  Since $f$ and $f^2$ are both
holomorphic on
$X-\cI(f)\cup f^{-1}(\cI(f))$, we see that $(f^2)^*T=(f^*)^2T$ on this set.  
Now $\cI(f^2)\subset
\cI(f)\cup f^{-1}(\cI(f))$, and let us write $V=\left(\cI(f)\cup
f^{-1}(\cI(f))\right)-\cI(f^2)$.  It suffices to show that $(f^2)^*T=0$ on $V$.  Since $T$
has codimension 1, it puts zero mass on any subvariety of codimension two.  Thus we
may suppose that $W$ is an irreducible component of $V$ of codimension 1.  Now by the
construction of $V$, we have $f(W-\cI(f))\subset\cI(f)$.  Thus $W$ is an exceptional
hypersurface.  Thus $\hat f(f(W-\cI(f)))$ is a subvariety of codimension at least
2.  It follows that $T$ puts no mass on this subvariety, and thus $(f^2)^*T$ puts
no mass on $W$.  We conclude, then that $(f^2)^*T$ puts no mass on $V$.  Since
$\Phi$ is obtained by extending by zero, we conclude that $\Phi_{f^2}T=\Phi_f^2T$. 
The proof for $n>2$ is similar. \ \ QED

Let us recall the cohomology group $H^{1,1}(X)$ which is given as the set of
smooth, $d$-closed  $(1,1)$-forms modulo $d$-exact 1-forms.  If $\omega$ is a
smooth $(1,1)$-form, then it acts on a $(d-1,d-1)$-form $\xi$ as
$\xi\mapsto\langle\omega,\xi\rangle=\int_{X}\omega\wedge\xi$.  Thus $\omega$
defines a (1,1)-current.  For $T\in\cT(X)$ there is a smooth (1,1)-form
$\omega_{T}$ such that the current $T-\omega_{T}$ is $d$-exact.   The
cohomology class $\omega_{T}\in H^{1,1}(X)$ is uniquely defined, so we have a
map $\cT(X)\to H^{1,1}(X)$.  If $V$ is a codimension 1 subvariety of $X$, we
let $\{V\}\in H^{1,1}(X)$ denote the cohomology class $\omega_{[V]}$
corresponding to the current of integration $[V]$.   The map $\Phi_{f}$ is
consistent with this passage to cohomology: $\{\Phi_{f}T\}=f^{*}\{T\}$.   We will say that $f$ is {\it
1-regular} if $(f^*)^n=(f^n)^*$ holds on $H^{1,1}(X)$. The
following is a consequence of Proposition 1.1:
\proclaim Proposition 1.2.  If (1.1) holds, then $f$ is 1-regular.

\bigskip
\centerline{\bf\S2. Degrees of the Iterates}
\medskip
\noindent A (complex) hyperplane $\cH\subset\P^d$, defined by a linear equation $\cH=\{\ell(x)=0\}$, gives a generator $H:=\{\cH\}$ of $H^{1,1}(\P^d)$.   If $V=\{h(x)=0\}\subset\P^{d}$ is a hypersurface of degree $m$, then $h/\ell^{m}$ is well-defined as a function on ${\bf P}^{d}$, so we have $[V]-[\cH]={1\over 2\pi}dd^{c}\log|h/\ell^{m}|$.  Thus $\{V\}=mH$, so the cohomology class of $\{V\}$ corresponds to the degree of $V$.  By definition, $\d(f)$ is the degree of the homogeneous polynomials defining $\hat f$, and we have $f^{*}H=\d(f)H$.  We will make use of the action on cohomology as a way of computing $\d(f)$.

The manifolds we will work with are obtained from ${\bf P}^{d}$ by blowing up points (the ``blowing up''  construction will be given in \S3).  This means that there is a sequence of spaces $X_{1}$, $X_{2}$,\dots,$X_{m}$ such that $X_{1}={\bf P}^{d}$ and $X_{m}=X$, and for each $j$ we have the following: there is a projection $\pi_{j}:X_{j}\to X_{j-1}$ and a finite set $S_{j-1}\subset X_{j-1}$ such that $\pi_{j}:X_{j}-\pi_{j}^{-1}S_{j-1}\to X_{j-1}-S_{j-1}$ is biholomorphic, and for each $s\in S$, the exceptional fiber is $\pi_{j-1}^{-1}s\cong{\bf P}^{d-1}$.  

For such $X$, we may describe $H^{1,1}(X)$ by the following inductive procedure.  
We start with $H$ as a basis of the (1,1)-cohomology of $X_{1}={\bf P}^{d}$.  Now suppose we have a basis $\cB_{j-1}$ for $H^{1,1}(X_{j-1})$.  We define a basis $\cB_{j}$ for $H^{1,1}(X_{j})$ by taking the elements $\pi_{j}^{*}b$ for $b\in\cB_{j-1}$, together with the classes of the exceptional fibers: $\{\pi_{j}^{-1}s\}$ for all $s\in S_{j-1}$.  

In particular, let us take $H_{X}:=\pi^{*}\{\cH\}$ as the first element of our basis $\cB$ of $H^{1,1}(X)$.  The rest of the basis elements can be taken to be exceptional fibers.  Let us suppose that the degree of $f$ is $m$.  Thus $f_{X}^{*}H_{X}= mH_{X} + E$, where $E$ denotes a sum over multiples of other basis elements from $\cB$.  The reason for the $m$ on the right hand side is as follows.  A generic line $L\subset{\bf P}^{d}$ does not intersect any of the centers of blow-up, so $\pi^{-1}$ is well defined in a neighborhood of $L$.  Thus $\pi^{-1}L$ intersects $f_{X}^{*}H_{X}$ with multiplicity $m$, since that is the multiplicity of intersection between $L$ and $f^{*}H$.

Now let $M$ be the matrix which represents $f_{X}^{*}$ with respect to the basis $\cB=\{H_{X},\dots\}$.  If $f_{X}$ is 1-regular, then $M^{n}$ is the matrix representation for $(f^{n}_{X})^{*}$ with respect to $\cB$.  Since $H_{X}$ is the first element of $\cB$, we have $d_n=(M^n)_{1,1}$.

There are various ways of representing $d_n$.  Let
$\lambda_1,\dots,\lambda_N$ denote the eigenvalues of the matrix $M$, and
let
$$\chi(x)=\prod_{j=1}^N(x-\lambda_j)=x^N+\chi_{N-1}x^{N-1}+\cdots+\chi_0$$
be the characteristic polynomial of $M$.   Let us first suppose that the
$\lambda_j$ are nonzero and have multiplicity one.  If we diagonalize $M$,
we find constants
$c_1,\dots,c_N$ such that
$$d_n=(M^n)_{1,1}=c_1\lambda_1^n+\cdots+c_N\lambda_N^n.\eqno(2.1)$$
Thus $\{d_n\}$ satisfies the recursion formula
$$d_{n+N}=\alpha_{N-1}d_{n+N-1}+\cdots+\alpha_0d_n\eqno(2.2)$$
where the coefficients $\alpha_j$ are determined by the characteristic
polynomial since we have $\alpha_j=-\chi_j$ for $0\le j\le N-1$.

Another way of producing the sequence $\{d_n\}$ is to find polynomials
$p(x)$ and $q(x)$ such that
$${p(x)\over q(x)}=\sum_{n=0}^\infty d_nx^n.$$
If we write $q(x)=\prod(x-r_j)$ and if $\d(p)<\d(q)$, then we may expand
$p(x)/q(x)$ into partial fractions we obtain
$${p(x)\over q(x)}=\sum_j {a_j\over r_j-x} =\sum_{n=0}^\infty
\left(\sum_j{a_j\over r_j^{n+1}}\right)x^n.\eqno(2.3)$$
Comparing (2.1) and (2.3), we see that after renumbering the $r_j$ if
necessary, we must have
$r_j^{-1}=\lambda_j$.  Thus $\chi$ and $q$ essentially determine each other:
$$q(x)\chi(0)=x^N\chi(1/x).\eqno(2.4)$$
Now let us suppose that the eigenvalues are nonzero (and not necessarily
simple).  Then we may approximate $M$ by a matrix $M'$ with simple
eigenvalues.  Equations (2.2) and (2.4) will hold for every such $M'$.  Since
the characteristic polynomial and the coefficients $((M')^n)_{1,1}$ depend
continuously on $M'$, equations (2.2) and (2.4) will continue to hold for
$M$.  It is not hard to adapt (2.4) to the case of zero eigenvalues.

\proclaim Theorem 2.1.  If $\delta>0$, then $\delta$ is the largest real zero of the characteristic polynomial $\chi(x)$.

\noindent{\it Proof.}  
Let $S$ denote the cyclic subspace spanned by $\{M^ne_1:n\ge0\}$, and let
$M_S$ be the restriction of $M$ to $S$.  For convenience, let us assume
first that $M_S$ is diagonalizable, and let $\lambda_1,\dots,\lambda_N$ be
the eigenvalues with corresponding eigenvectors $v_1,\dots,v_N$.  We may
write $e_1=\sum c_jv_j$ and $M^ne_1=\sum c_j\lambda_j^nv_j$.  Since $e_1$
is cyclic on $S$, there are nonzero numbers $a_j$ such that
$$d_n=e_1\cdot M^ne_1=\sum_{j=1}^Na_j\lambda^n_j.\eqno(2.5)$$
Now at least one of the $\lambda_j$ must have modulus $\delta$ (the
spectral radius of $M_S$).  We claim that since $d_n\ge0$, it follows that
$\delta$ itself must be an eigenvalue.

Let us suppose, by way of contradiction, that $\delta\ne\lambda_j$ is not
an eigenvalue.  It is easy to reduce to the case where all the eigenvalues
in (2.5) have modulus $\delta$.  Let us write $\lambda_j=\delta e^{2\pi
i\theta_j}$ with $0<\theta_j<1$.  We consider two cases separately.  The
first case is that all $\theta_j$ are rational.  Thus all the $\lambda_j$
are $M$th roots of unity for some $M$.  It follows that
$\sum_{n=1}^M\lambda_j^n=0$, and $d_n\ge0$, so we must have $d_n=0$ for all
$n$.  However, since the vectors $\{(\lambda_1^n,\dots,\lambda_N^n),1\le
n\le N\}$ are a basis for ${\bf C}^N$, the condition $d_n=0$ for $1\le n\le
M$ implies that the
$a_j$ all vanish, which is a contradiction.

The other case is where the $\theta_j$ are  not all rational.  We consider
the closed (Lie) subgroup $G$ of ${\bf R}^N/{\bf Z}^N$ generated by
$\{n(\theta_1,\dots\theta_N):n\in{\bf Z}\}$.  Let us define
$h(\xi_1,\dots,\xi_N)=\sum a_j e^{2\pi i \xi_j}$ for
$(\xi_1,\dots,\xi_N)\in{\bf R}^N/{\bf Z}^N$.  By continuity, we have
$h(\xi)\ge0$ for $\xi\in G$. Since
$G$ is a closed subgroup, it contains a point $(\phi_1,\dots,\phi_N)$ with
rational coordinates.  Thus $h(n\xi)\ge 0$ for all $n\ge0$.  Arguing as
before, we reach the same contradiction that the $a_j$ all vanish.

Finally, if $M_S$ is not diagonalizable, we may decompose it into cyclic
subspaces corresponding to the different eigenvalues, and the $a_j$ in the
formula above become polynomials of $n$.  Now we may replace the
polynomials by the highest degree coefficients and proceed as before.
\bigskip
\centerline{\bf \S3.  1-Regularization of $J$}
\medskip
\noindent A basic map we work with is defined by
$$ J[x_0 :x_1:\cdots : x_d ] = [ {x_0}^{-1}:{x_1}^{-1}: \cdots
:{x_d}^{-1}] = [x_{\hat 0}:x_{\hat 1}:\cdots :x_{\hat d}] $$ 
where we write $ x_{\hat \jmath} = \prod _{i \neq j} x_j.$  This is an involution because
$$J^{2}[x_{0}:\cdots:x_{d}] = [(x_{0}\dots x_{d})^{d-1}x_{0}:\cdots:(x_{0}\dots x_{d})^{d-1}x_{d}] =[x_{0}:\cdots:x_{d}]$$
for $x_{0}\dots x_{d}\ne 0$.  In order to discuss the behavior of $J$, we introduce some notation.
For a subset $ I \subset \{ 0,1, \dots , d \}$, we write the complement
as $\hat I=\{0,\dots,d\}-I$.  Let us set 
$$
\Sigma_I  = \{x\in\P^d: x_i = 0, \ \forall i \in I \},\ \  \Sigma^*_I = \{
x\in\Sigma_I: x_j\neq 0, \  \forall j \in\hat I \}.
$$
Thus $\Sigma_I$ is closed, and $\Sigma_I^*$ is a dense, open subset of $\Sigma_I$.
We see that the indeterminacy locus is given by $\cI(J)=\bigcup_{|I|\ge2}\Sigma_I$.  In fact, we have a stratification given by
$$ {\cI}(J) = \bigcup ^d _{j=0} ( \Sigma_j - \Sigma_j^*)
=\bigcup_{|I|\ge2}\Sigma_I^*.$$  
The action of $J$ corresponds to the involution $I\leftrightarrow \hat I$ of the set of subsets of $\{0,\dots,d\}$: for $|I|\ge2$, $\hat J$ acts as
$$\hat J:\Sigma^*_I\ni p\mapsto \hat J(p)=\Sigma_{\hat I}.$$  

The points of ${\bf P}^d-\cI(J)$ where $J$ is not a local diffeomorphism is $\bigcup_{j=0}^d\Sigma_j^*$.  We use the notation
$$e_j:=\Sigma_{\hat \jmath}=\{[0:\cdots:0:1:0:\cdots:0]\},$$
so $J(\Sigma^{*}_{j})=e_{j}$, and the exceptional locus is $\cE=\Sigma_{0}\cup\dots\cup\Sigma_{d}$.

$J$ is not 1-regular, since we have
$\hat J:\Sigma_j^*\to e_j\to\Sigma_j$.  We show here how a 1-regularization of $J$ may be obtained by
blowing up.  Let $\pi:X\to\C^{d}$ denote the space
$\C^d$ blown up at the origin.  We represent the blow-up as
$$X=\{((z_1,\dots,z_d),[\xi_1:\cdots:\xi_d])\in\C^d\times\P^{d-1}:
z_i\xi_j=z_j\xi_i, 1\le i,j\le d\},$$
and $\pi(z,\xi)=z$.
Thus $X$ is a smooth $d$-dimensional submanifold of
$\C^d\times\P^{d-1}$.  Let $E:=\pi^{-1}(0)$ denote the fiber over the
origin.  Thus $E\cong\P^{d-1}$, and $\pi:X-E\to\C^d-\{0\}$ is
biholomorphic; the inverse map is given by $(z_1,\dots,z_d)\mapsto ((z_1,\dots,z_d),[z_1:\cdots:z_d])$, for $z\ne(0,\dots,0)$.  If $V$ is a complex subvariety of ${\bf C}^d$, we identify it as a subvariety $V_{X}\subset X$ as follows: by $V_{X}$, we mean the closure of $\pi^{-1}(V-\{0\})$ inside $X$.  Thus if $0\in V$, this is a proper subset of $\pi^{-1}V=V_{X}\cup E$.  When there is no danger of confusion, we will write $V$ for $V_{X}$.

Let us see how the operation of blow-up modifies the map $J$.  We may identify a
neighborhood of
$e_0=[1:0:\cdots:0]=J(\Sigma^*_0)$ inside
$\P^d$ with
$\C^d$ via the map $ [1:z_1:\cdots:z_d]\leftrightarrow(z_1,\dots,z_d)$.  Performing the blow-up 
$\pi:X\to\C^d$ at $e_{0}$ in the range of $J$ induces a (partially defined) map $ J_{X}:{\bf C}^d\to X$, given by
$$J_{X}:=\pi^{-1}\circ J:\{x\in\P^d:x_1,\dots,x_d\ne0\}\to X.$$  
For $x_0\ne0$, we have
$$J_{X}:[x_0:x_1:\cdots:x_d]\mapsto[x_0^{-1}:x_1^{-1}:\cdots:x_d^{-1}]=
[1:{x_0\over x_1}:\cdots:{x_0\over x_d}]$$
$$\leftrightarrow 
\left(z_1={x_0\over x_1},\dots,z_d={x_0\over x_d}\right) \mapsto
\pi^{-1}(z_1,\dots,z_d)\in X.$$
Thus, letting $x_0\to0$, we obtain the map $J_{X}:\Sigma^*_0\to
E\cong\P^{d-1}$ given by
$$J_{X}(0:x_1:\cdots:x_d)=[\xi_1=x_1^{-1}:\cdots:\xi_d=x_d^{-1}].$$ We
see that $J_{X}$ is a local diffeomorphism at points of $\Sigma^*_0$.

 The process of blowing up a point is in fact local and can
be performed at any point of a complex manifold.  Let $\pi:X\to{\bf P}^d$ denote the complex manifold obtained by blowing up at the centers $\{e_0,\dots,e_d\}$, and let $E_{j}=\pi^{-1}e_{j}$ denote the exceptional fiber over $e_{j}$.   
 
 Let us describe the induced birational map  $J_X:X\to X$.  Since $$\pi:X-\bigcup_{j=0}^{d}E_{j}\to {\bf
P}^{d}-\{e_{0},\dots,e_{d}\}$$ is a biholomorphism, it follows that $\cI(J_{X})\cap (X-\bigcup
E_{j})=\bigcup_{|I|\ge2}\Sigma_{I}\cap(X-\bigcup E_{j})$.  Further, the calculation above showed that 
$J_{X}|\Sigma^{*}_{j}$ is essentially $J$, and thus $J_{X}|E_{j}$ is essentially $J$ on ${\bf P}^{d-1}$. 
Thus we conclude that $\cI(J_{X})=\bigcup_{|I|\ge2}\Sigma_{I}$, where $\Sigma_{I}\subset X$ is
interpreted as above.  In particular, $\cI\cap E_{j}$ has codimension 2 in $E_{j}$.  Now the restriction
of $J_{X}$ to $X-\bigcup\Sigma_{j}$ may be identified with the restriction of $J$ to ${\bf
P}^{d}-\bigcup\Sigma_{j}$, which is a diffeomorphism.   We have also seen that $J_{X}$ is a local
diffeomorphism on $\bigcup\Sigma^{*}_{j}$.  Thus $J_{X}$ is a local diffeomorphism at all points of
$X-\cI(J_{X})$.  This means that the exceptional locus is empty, and thus $J_X$ is 1-regular.   

If $d=2$,
then
$J_{X}$ is in fact holomorphic, i.e., $\cI(J_{X})=\emptyset$.  If $2\le |I|<d$, then $\Sigma_{I}\cap({\bf
P}^{d}-\{e_{0},\dots,e_{d}\})\ne\emptyset$, so $\Sigma_{I}$ intersects the exceptional fibers $E_{j}$ for
all $j$ such that $j\notin I$.

Let us remark that a linear map $L$ also induces a birational map
$L_X=\pi^{-1}\circ L\circ\pi:X\to X$.  It is evident that
$\cI(L_X)\subset\{e_0,\dots,e_d\}$ and $\cE(L_X)\subset E_0\cup\dots\cup E_d$.  If
$L^{-1}e_j=e_i$ for some $i$ and $j$, then $L_X$ is biholomorphic in a neighborhood
of $E_i$ and maps $E_i$ to $E_j$.  If $L^{-1}e_j$ is not one of these points $e_i$,
then $L^{-1}e_j\in\cI(L_X)$, and $\hat L(L^{-1}e_j)=E_j$.  And if $Le_j$ is not one
of the $e_i$, then $L_XE_j$ is a point, so $E_j\subset\cE(L_X)$.
We see that $L_X$ fails to be 1-regular exactly when there is a point $e_i$ such
that $Le_i$ is not one of the $e_j$'s, but $L^ne_i=e_j$ for some $n\ge2$ and some
$j$.

From the discussion of the previous paragraph, we can deduce the action of $L^*_X$
on $H^{1,1}(X)$.  Namely, let $H_X$ denote the cohomology class of a
hyperplane.  Since neither $H$ and $LH$ will contain any of the $e_i$'s for generic
$H$, we see that $L_X^*H_X=H_X$.  Further, we have $L_X^*E_i=E_j$ for the pairs
$(i,j)$ such that
$Le_j=e_i$.  Let $M$ be the $(d+1)\times(d+1)$ matrix such that $m_{i,j}=1$ if $Le_j=e_i$ and 0
otherwise.  Then, with respect to the basis
$\cB= \{H_X,E_0,\dots,E_d\}$ of $H^{1,1}(X)$, we have
$$L_X^*=\pmatrix{1&0\cr
0&M\cr}.$$

Next we discuss the action induced by $J_{X}$ on $H^{1,1}(X)$.  Let $H \in H^{1,1}({\bf P}^{d})$ denote the class of a hyperplane, and let $H_{X}=\pi^{*}H$ denote the induced class in $H^{1,1}(X)$.  When there is no danger of confusion, we will also denote $H_{X}$ simply by $H$.  Let $E_{j}\in H^{1,1}(X)$ denote the cohomology class induced by $E_{j}$.  Thus $\{H_{X},E_{0},\dots,E_{d}\}$ is a basis for $H^{1,1}(X)$, and we will represent $J^{*}_{X}$ as a matrix with respect to this basis.  

Let $\{\Sigma_{0}\}\in H^{1,1}(X)$ denote the class induced by $\Sigma_{0}$.  We wish to represent $\{\Sigma_{0}\}$ in terms of our basis.  Let us start by observing that $\Sigma_{0}$ is a hyperplane in ${\bf P}^{d}$, and so $\Sigma_{0}=H\in H^{1,1}({\bf P}^{d})$.  Thus we have
$H_{X}=\pi^{*}H=\pi^{*}\Sigma_{0}$.  By our formula for the pullback of a current, we have that $\pi^{*}\Sigma$ will correspond to the current of integration over $\pi^{-1}\Sigma_{0}=\Sigma_{0}\cup E_{1}\cup\dots\cup E_{d}$.  (Since $e_{0}\notin\Sigma_{0}$, the divisor $E_{0}$ will not be involved.)  It remains to determine the multiplicities of the different components.  Now on the set $x_{i}\ne0$, the current $\{\Sigma_{0}\}$ is represented by the potential $\log|x_{0}/x_{i}|$.  Let us choose $i=d$ for convenience, and use affine coordinates $z_{0}=x_{0}/x_{d},\dots,z_{d-1}=x_{d-1}/x_{d}$.  In these coordinates, the potential for $\{\Sigma_{0}\}$ is given by $h:=\log|z_{0}|$.  Let us work in a neighborhood of the exceptional fiber $E_{k}$ for some $1\le k\le d-1$.   On the dense open subset $\xi_{0}\ne0$, we may write a point of the fiber as $[1:\xi_{1}:\dots:\xi_{d-1}]$.  With this, we define an affine coordinate system 
$$(z_{0},\xi_{1},\dots,\xi_{d-1})\mapsto ((z_{0},\xi_{1}z_{0},\dots,1+\xi_{k}z_{0},\xi_{k+1}z_{0},\dots,\xi_{d-1}z_{0}),[1:\xi_{1}:\dots:\xi_{d-1}]).$$  
In this coordinate system, we see that $E_{k}$ is given by $z_{0}=0$.  The potential for the current $\pi^{*}\{\Sigma_{0}\}$ is then given by
$\pi^{*}h = h\circ\pi$.  It follows that the multiplicities are one, so
$$\pi^{*}\Sigma_{0}= \{\Sigma_{0}\} + \sum_{i\ne0}E_{i}.$$
Combining this with the previous equation, we obtain
$$\{\Sigma_{j}\}=H_{X}-\sum_{i\ne j}E_{i}.$$

We have seen that $J_{X}$ is a diffeomorphism from $\Sigma_{j}^{*}$ to its image in $E_{j}$.  Since
$J_{X}$ induces a diffeomorphism outside a subvariety of codimension 2, and we are pulling back
cohomology classes of codimension one, it follows that $$J_{X}^{*}E_{j}=\{\Sigma_{j}\}=H_{X}-\sum_{i\ne
j}E_i.$$

Next we need to determine $J_{X}^{*}H_{X}$.  A generic hyperplane ${\cal H}$ in ${\bf P}^{d}$ does not meet any of the $e_{j}$ and may be considered to be a subset of $X$.  Thus it generates $H_{X}$.  Thus we  consider the restriction $J|_{X-\cI}$ and determine the class $\{(J^{-1}{\cal H})-\cI\}=J^{*}_{X}H_{X}\in H^{1,1}(X)$.  Let us start with the observation which connects $H$, $H_{X}$, and the preimage of  ${\cal H}$:
$$d\cdot H_{X} = \pi^{*}(d\cdot H) = \pi^{*}(J^{*}H) = \pi^{*}\{J^{-1}{\cal H}\}.$$
A hyperplane has the form ${\cal H} = \{h=0\}$ for some $h=\sum a_{j}x_{j}$.  Thus $J^{-1}{\cal H} = \{\sum a_{j}x_{\hat\jmath}=0\}$, and  $\log h\circ J\circ\pi$ will be a potential for $\pi^{*}\{J^{-1}{\cal H}\}$.  The element $ \pi^{*}\{J^{-1}{\cal H}\}\in H^{1,1}(X)$ will be $\{J^{-1}{\cal H}\}$ and a linear combination of the $E_{j}$.  We need only determine the multiplicities of the $E_{j}$.  Let us consider $E_{d}$.  Since $x_{d}\ne0$, we work in an affine coordinate system $(z_{0},\dots z_{d-1})$. We write points in the fiber as $[\xi_{0}:\dots:\xi_{d-1}]$.  On a dense open subset of the fiber we  have $\xi_{0}\ne0$, and so with the same coordinate system as above we have
$$h\circ J\circ \pi =A(\xi)z_{0}^{d-1}+B(\xi)z_{0}^{d}$$
For generic $\xi$, $A(\xi)\ne0$, so this vanishes to order $d-1$ in $z_{0}$, and thus the multiplicity of $E_{d}$ (and all $E_{j}$) is $d-1$.  This gives
$$\pi^{*}\{J^{-1}{\cal H}\}=\{J^{-1}{\cal H}\} + \sum_{j}(d-1)E_{j}.$$
Finally, we use the fact that $J_{X}^{*}H_{X}=\{J^{*}{\cal H}-\cI\}$
to conclude that
$$J_{X}^{*}H_{X} = d\cdot H_{X}+\sum_{j}(1-d)E_{j}.$$
Thus we may write the action on $H^{1,1}$ with respect to our basis in matrix form: 
$$\pmatrix{  d & 1 &1& \dots & 1 \cr
1-d& 0 &-1 & \dots & -1 \cr
1-d&-1&0 &\dots& -1\cr
\vdots&\vdots&\vdots&\ddots& \vdots\cr
1-d&-1&-1& \dots &0\cr}\eqno(3.1)$$
The fact that this matrix is an involution corresponds to the fact that $J_{X}$ is 1-regular.
\bigskip
\centerline{\bf \S4. Elementary Mappings}
\medskip\noindent
This Section is devoted to a discussion of mappings of the form $f=L\circ J$, where $L$ is a linear map of ${\bf P}^{d}$, and $J$ is as in the previous section.  For $p\in X$, we define the orbit $\cO(p)$ as follows.  If $p\in\cE\cup\cI$, then $\cO(p)=\{p\}$.  If there is an $N\ge1$ such that $f^{j}p\notin\cE\cup\cI$ for $0\le j\le N-1$ and $f^{N}p\in\cE\cup\cI$, then $\cO(p)=\{p,f p,\dots,f^{N}p\}$.  Otherwise, we have $f^{j}p\notin\cE\cup\cI$ for all $j\ge0$, and we set $\cO(p)=\{p,fp,f^{2}p,\dots\}$.  In the first two cases (when the orbit is finite), we say that the orbit is {\it singular}.  Otherwise, we say that the orbit is {\it nonsingular}.

The orbits starting at the image points of the exceptional hypersurfaces $\Sigma_{j}$, for $0\le j\le d$ have special importance.  We will use the notation $\alpha_{j}=f\Sigma^{*}_{j}=Le_{j}$ for the image of the $j$th exceptional hypersurface, which is identified with $j$th column of the matrix $L$; and we let $\cO_{j}:=\cO(f\Sigma^{*}_{j})=\cO(\alpha_{j})$ denote its orbit.
We say that the mapping $f$ is {\it elementary} if for each $0\le j\le d$, the orbit $\cO_{j}$  is either nonsingular or, if it is singular, it ends at one of the points $\{e_{0},\dots,e_{d}\}$.  

Now suppose that $f$ is elementary.  We define an {\it orbit list} to be a list of singular orbits of exceptional components $\cO_{i}=\cO(\alpha_{i})$ with sequential indices:
$$\cL=\{\cO_{a},\cO_{a+1},\cO_{a+2},\dots,\cO_{a+\mu}\}$$
such that if $0\le j<\mu$, then the endpoint of $\cO_{a+j}$ is $e_{a+j+1}$.  In other words, if $j<\mu$,
then the orbit
$\cO_{a+j}$ must be singular, and the ending index $k$ of the endpoint $e_{k}$ of this orbit is the beginning
index  of the next orbit in the list.  Let us suppose that the last orbit in the list, $\cO_{a+\mu}$, ends at the point $e_k$.  We say that the list $\cL$ is {\it open} if the orbit $\cO_k=\cO(\alpha_k)$ is nonsingular.  We say that $\cL$ is {\it closed}  if $k=a$.

Renumbering the variables, if necessary,  we may group the orbits into maximal orbit lists $\cL_{1},\dots,\cL_{\nu}$.  It follows from the maximality that each $\cL_{j}$ is either open or closed.  Let us define $\cA$ to be the set of indices $i$ such that $\cO_{i}$ is a singular orbit and is the first orbit in an open orbit list.  Let $\Omega$ consist of the indices $j$ such that $e_{j}$ is the endpoint of a singular orbit.

Now we construct the 1-regularization of $f$.  Let $S=\{i:\cO_i{\rm\ is\ singular}\}$, and let
$\cO_S:=\bigcup_{i\in S}\cO_i$.  Let
$\pi:X\to{\bf P}^{d}$ be the space obtained by blowing up each of the points of $\cO_S$.  For
$p\in\cO_S$, we let $\cF_{p}$ denote the exceptional fiber $\pi^{-1}p$ in $X$
over $p$, and we also let $\cF_{p}$ denote the induced cohomology class in $H^{1,1}(X)$.  Repeating the
reasoning of the previous section, we see that the hypersurfaces $\Sigma_{j}\subset X$, $j\in S$,
are not exceptional for the induced birational map $f_{X}:X\to X$.  Thus $f_{X}$ is 1-regular.

Let us determine the induced mapping $f_{X}^{*}$ on $H^{1,1}(X)$.  The class $H_{X}$, together with the classes
$\cF(p)$ for $p\in\cO_S$, form a basis for $H^{1,1}(X)$.   For $i\in S$, we have
$$\Sigma_{i}\to f\Sigma^{*}_{i}=\alpha_{j}\to\cdots\to f^{n_{j}-1}\alpha_{j}= f^{n_{i}}\Sigma^{*}_{i}=e_{\beta_{i}}.$$
At each of the points $f^{j}\alpha_{i}$, $0\le j\le n_{i}-2$, $f$ is locally biholomorphic, so $f_{X}$ induces a biholomorphic map of a neighborhood of the fiber $f_{X}:\cF_{f^{j}\alpha_{i}}\to \cF_{f^{j+1}\alpha_{i}}$.  We conclude that 
$$f^{*}_{X}\cF_{f^{j+1}\alpha_{i}}=\cF_{f^{j}\alpha_{i}}{\ \ \rm for\ }0\le j\le n_{i}-1\eqno(4.1)$$
and
$$f^{*}_{X}\cF_{\alpha_{i}}=\{\Sigma_{i}\}$$
where  $\{\Sigma_{i}\}$ denotes the class induced by $\Sigma_{i}$ in $H^{1,1}(X)$.  As in the previous section, we have
$$\{\Sigma_{i}\}=H_{X}-\sum_{p} \cF_{p},$$
where the sum is taken over all blow-up centers which belong to $\Sigma_{i} \cap \cI$.  The set of blow-up centers which belong to $\cI$ is $\Omega$, and the only question is whether $i\in\Omega$.
In fact, we have $i\in\Omega$ if $i\notin\cA$.  Thus we have
$$\eqalign{ f_{X}^{*}\cF_{\alpha_{i}}=& H_{X} -\cF_{\Omega}+\cF_{e_{i}}{\rm\ \ \ \ \ \ }i\notin\cA,\cr 
f_{X}^{*}\cF_{\alpha_{i}}= &H_{X} -\cF_{\Omega}{\rm\ \ \ \ \ \ \ \ \ \ \ \ \ \  }i\in\cA,\cr}\eqno(4.2)$$
where we have adopted the notation $\cF_{\Omega}:=\sum_{t\in\Omega}\cF_{e_{t}}$.
Finally, to pull back the class of a hyperplane, we use the fact that $L$ is biholomorphic, and for a generic hyperplane $\cH$, the preimage $L^{-1}\cH$ is again a generic hyperplane.  Thus we may use $f_{X}^{-1}=J_{X}^{-1}\circ L^{-1}$ and argue as in the previous section to find:
$$f^{*}_{X}H_{X}=\{f_{X}^{-1}\cH\} =\{J_{X}^{-1}\cH\}=dH_{X}
+(1-d)\cF_{\Omega}.\eqno(4.3)$$

Let $\cL^{c}=\{\{\cO_{a_{1}},\dots,\cO_{a_{1}+\mu_{1}}\},\{\cO_{a_{2}},\dots,\cO_{b_{2}+\mu_{2}}\},\dots,\{\cO_{a_{m}},\dots,\cO_{a_{m}+\mu_{m}}\}\}$ denote a listing of the set of closed orbit lists, and let $\cL^{o}$ be a listing of the set of those open orbit lists which contain singular orbits.   For an orbit  $\cO$, we let $|\cO|$ denote its length, and by $\#\cL^{c}$ we denote the set of lists of lists of orbit lengths $$\#\cL^{c}=\{\{|\cO_{a_{1}}|,\dots,|\cO_{a_{1}+\mu_{1}}|\},\dots,\{|\cO_{a_{m}}|,\dots,|\cO_{a_{m}+\mu_{m}}|\}\}.$$  
We see that the mapping $f^{*}_{X}$ is determined by $\#\cL^{c}$ and $\#\cL^{0}$.  Thus we have the following:
\proclaim Theorem 4.1.  If $f=L\circ J$ is elementary, then the dynamic degree $\delta(f)$ is determined by
$\#\cL^{c}$ and $\#\cL^{o}$.

Henceforth, we will abuse notation and simply write $\cL^c$ and $\cL^o$ for the orbit list structure $\#\cL^{c}, \#\cL^{o}$, since only the lengths of the orbits
(and not the specific points) are used in computing $f^*$ and
$\delta$.  Thus we may consider
$\delta(\cL^c,\cL^o)$ to be a number which is determined by two
sets of lists of positive integers.  In addition, the
characteristic polynomial may be explicitly computed in terms of
$\cL^c$ and $\cL^o$.  This is done in the Appendix.  By Theorem
2.1,  $\delta$ is the largest real zero of $\chi(x)$, and by
Theorem A.1 $\chi(d)>0$, so we have the following:
\proclaim Theorem 4.2.   If $f$ is an elementary mapping of ${\bf P}^d$ with at least one singular orbit, then $\delta<d$.

\bigskip
\centerline{\bf\S5.  Comparison Results}
\medskip\noindent  We saw in the previous section that for an elementary mapping, $\delta$ is determined by the orbit list structure $\cL^{c},\cL^{o}$.  Here we develop some results which may be interpreted as giving monotonicity properties of this dependence, or equally well, as giving a method of comparing $\delta$ whenever the orbit lists may be compared.   Let us describe how to compare orbit lists and lists of lists.  Our first comparison theorem involves lists with the same structure pattern but different orbit lengths.  Let $\cL=\{N_{1},\dots,N_{\ell}\}$ and $\hat\cL=\{\hat N_{1},\dots,\hat N_{\hat\ell}\}$  be the structures of two orbit lists.  We say that $\hat\cL$ {\it has longer orbits than} $\cL$ if $\ell=\hat\ell $ and $|\hat N_{i}|\ge|N_{i}|$ for all $1\le i\le \ell$.  If these are closed orbit lists, we may also allow circular permutations of the orbits in our comparison.  Now if $\cL=\{\cL_{1},\dots,\cL_{\mu}\}$ and  $\hat\cL=\{\hat\cL_{1},\dots,\hat\cL_{\hat\mu}\}$ 
are lists of lists, we say that $\hat\cL$ {\it has longer orbits than} $\cL$ if they are both of the same type (either open or closed) and if $\hat\mu=\mu$, and after a possible permutation of the index set $\{1,\dots,\mu\}$, the list  $\hat\cL_{j}$ has longer orbits than $\cL_{j}$ for each $1\le j\le\mu$.

\proclaim Theorem 5.1.  Let $f$ and $\hat f$ be two elementary maps of ${\bf P}^{d}$.  Let $\cL^{o},\cL^{c}$ (respectively, $\hat\cL^{o},\hat\cL^{c}$) be the orbit list structure of $f$ (respectively, $\hat f$).  If $\hat\cL^{c},\hat\cL^{o}$ has longer orbits than $\cL^{c},\cL^{o}$, then $\delta(\hat f)\ge \delta(f)$.  If $\delta(f)>1$, then the inequality is strict.

\noindent{\it Proof.}   Proceeding by induction, we may assume that all the orbit lengths except one are the same.  We will suppose that the orbit length changes inside one of the closed orbit lists.  (The proof of the case if the orbit list is open is similar.)  Without loss of generality, we may suppose that the orbit which is changed is the first orbit inside $\cL^c_1$, and its length is $N_{1,1}$, and that the orbit length is $\hat N_{1,1}=N_{1,1}+1$ inside $\hat\cL_1^c$.  Let $\chi$ (respectively $\hat\chi$) denote the characteristic polynomial corresponding to $\cL^c,\cL^o$ (respectively $\hat\cL^c,\hat\cL^o$).  Now recall from (A.3) that the characteristic polynomial has the general form
$$\chi(x)=(x-d) T_1{\prod}'T+(x-1)S_1{\prod}' T+(x-1)T_1\sum S({\prod}''T) .\eqno(5.1)$$
Since the orbit lists agree except at the first orbit of the first list, we have $T_i=\hat T_i$ and $S_i=\hat S_i$ except for $i=1$.  For $x>1$, we have
$$\gamma:={\hat T^c_1(x)\over T^c_1(x)} = {x\prod_{j=1}^{\ell^c_i}x^{N^c_{i,j}}-1\over \prod_{j=1}^{\ell^c_i}x^{N^c_{i,j}}-1} = x + {x-1\over \prod_{j=1}^{\ell^c_i}x^{N^c_{i,j}}-1}>x.$$
Similarly, we find that for $x>1$, we have
$$\rho:= {\hat S^c_1(x)\over S^c_1(x)} = x- {(x-1)\cdot({\rm positive \ terms})\over S^c_1(x)} < x.$$
It follows that $\rho<\gamma$.  Substituting $\rho$ and $\gamma$ into (5.1), we find that if $1<x<d$, then
$$\hat\chi(x) = \gamma(x-1)T_1{\prod}'T + \rho(x-1)S_1{\prod}'T+\gamma(x-1)T_1\sum S{\prod}''T$$
and thus $\hat\chi(x)<\gamma\chi(x)$, since $T_i,S_i>0$, and $\rho<\gamma$.

Finally, if we set $x=\delta$, then by Theorem 2.1, we have $\chi(\delta)=0$, which gives $\hat\chi(\delta)<0$.  Thus the largest root of $\hat\chi$ will be greater than $\delta$.  This gives us the desired result.  \ \ QED\medskip

Next we discuss the limiting behavior as the length of one (or several) of the orbits becomes unbounded.

\proclaim Theorem 5.2.  Let $f$ be an elementary map of ${\bf P}^{d}$ with orbit structure $\cL^{c},\cL^{o}$ and with $\delta(f)>1$.  Let $\cL=\{N_{1},\dots,N_{\ell}\}$ be one of the orbit lists in this structure, and let
$$\delta_{i}:=\lim_{N_{i}\to\infty}\delta(f).$$
Then $\delta_{i}$ is the dynamical degree corresponding to the orbit list structure $\hat\cL^{c},\hat\cL^{o}$ which is obtained as follows:
\vskip0pt{Closed case:}  If $\cL$ is a closed orbit list, then $\hat\cL^{c}$ is obtained by deleting $\cL$ from $\cL^{c}$ and adding the open list $\{N_{i+1},\dots,N_{\ell},N_{1},\dots,N_{i-1}\}$ to $\cL^{o}$.
\vskip0pt{Open case:}  If $\cL$ is an open orbit list, then we set $\hat\cL^{c}=\cL^{c}$ and replace $\cL$ in $\cL^{o}$ according to the following cases:
\vskip0pt\noindent If $i=1$, replace $\cL$ by the open list $\{N_{2},\dots,N_{\ell}\}$.
\vskip0pt\noindent If $1<i<\ell$, replace $\cL$ by the pair of open lists $\{N_{1},\dots,N_{i-1}\}$ and $\{N_{i+1},\dots,N_{\ell}\}$.
\vskip0pt\noindent If $i=\ell$, replace $\cL$ by the open list $\{N_{1},\dots,N_{\ell-1}\}$.

\noindent{\it Proof.}  There are four cases to consider.  The proofs of all these cases are similar, so we consider only the first case.  Since $\cL$ is closed, we may perform a circular permutation so that we have $i=\ell$.   Let $\chi(x)$ denote the characteristic polynomial as given by the formula (A.3), and let $\hat\chi$ denote the characteristic polynomial for the orbit structure obtained from $\cL^c,\cL^o$ by replacing the list $\{N_1,\dots,N_\ell\}$ by $\{N_1,\dots,N_{\ell-1}\}$.   Inspecting the formula (A.3), we may write
$$x^{-N_\ell}\chi(x)= \hat\chi(x) + O(x^{-N_\ell})\eqno(5.2)$$
for $x>1$.  For each value of $N_\ell$, we let $\delta_{N_\ell}$ denote the corresponding dynamical degree, which is also the largest real zero of $\chi$.  By Theorem 5.1, $\delta_{N_\ell}$ is monotone increasing.  Thus $\delta_{N_\ell}^{-N_\ell}\to0$.  We conclude that the $O$ term in (5.2) vanishes as $N_\ell\to\infty$, and so the limiting value, $\delta_i$ is the largest real zero of $\hat\chi$. \ \ QED\medskip

\proclaim Theorem 5.3.  Let $\cL^c, \cL^o$ be the orbit list structure of an elementary map of ${\bf P}^d$.  If we let $\hat\cL^c,\hat\cL^o$ be the orbit list structure obtained by adding an orbit list to $\cL^c$ or $\cL^o$, then $\delta(\hat\cL^c,\hat\cL^o)\le\delta(\cL^c,\cL^o)$.

\noindent{\it Proof.}  For the new orbit list, let $T(x)$ and $S(x)$ denote the polynomials corresponding to the definitions in (A.1-2).  Let $\chi$ denote the characteristic polynomial corresponding to the old orbit list structure, and let $\hat\chi$ denote the characteristic polynomial corresponding to the new one.  Thus we have
$$\eqalign{\hat\chi(x) &= (x-d)T(x){\prod}'T + (x-1)S(x){\prod}'T+(x-1)T(x)\sum S{\prod}''T\cr
&=T(x)\chi(x)+(x-1)S(x){\prod}'T,\cr}$$
where the notation ${\prod}'$ means we are taking the product over all of the polynomials $T_i^c$ and $T_i^o$, except the new $T(x)$.  If we let $x=\hat\delta$ be the largest zero of $\hat\chi$, then we have
$$0=T(\hat\delta)\chi(\hat\delta)+(\hat\delta-1)S(\hat\delta)[\cdots].$$
Since $\hat\delta\ge1$, it follows that all the terms except $\chi(\hat\delta)$ on the right hand side of the equation are positive, so $\chi(\hat\delta)\le 0$.  Thus the largest zero of $\chi$ is greater than or equal to $\hat\delta$. \ \ QED \medskip

\proclaim Theorem 5.4.  Let $\cL^c,\cL^o$ denote the orbit list structure of an elementary map of ${\bf P}^d$, and let $\cL=\{N_1,\dots,N_\ell\}$ denote the structure of one of the lists.  For $1\le j<\ell$ there is a number $M^*=M^*(j,\cL,\cL^c,\cL^o)$ with the following property: Given $M$, we let $\cL(j)=\{N_1,\dots,N_j,M,N_{j+1},\dots,N_\ell\}$ denote the list obtained by adding an orbit of length $M$ at the $(j+1)$st place in the list $\cL$.  Let $\hat\cL^c,\hat\cL^o$ denote the new orbit list structure obtained by replacing $\cL$ with $\cL(j)$.  Then if $M<M^*$, we have $\delta(\hat\cL^c,\hat\cL^o)<\delta(\cL^c,\cL^o)$; if $M>M^*$, we have 
$\delta(\hat\cL^c,\hat\cL^o)>\delta(\cL^c,\cL^o)$.

\noindent{\it Remark.}  In the open case, if $j=0$ or $j=\ell$, then by Theorems 5.1 and 5.2 we can only reduce $\delta$ by adding an orbit in the $j$th position; this means that $M^*=\infty$.

\noindent{\it Proof.}  Let us assume that $\cL$ is an open orbit list.  (The proof for the case of a closed orbit list is similar.)  Without loss of generality we may suppose that the orbit list $\cL$ is the orbit list inside $\cL^o$.  Let us define 
$$\eqalign{\varphi(x,M)= x^M[1 + & \sum_{i=1}^{j-1} \prod_{k=1}^i x^{N_k}] [1 + \sum_{i=j+2}^\ell\prod_{k=i}^\ell x^{N_k}]\cr
& - [1+\sum_{i=1}^j \prod_{k=1}^i x^{N_k}] [1 + \sum_{i=j+1}^\ell \prod_{k=i}^\ell x^{N_k}].}$$
For the new orbit list structure we have 
$$\hat S_1(x) = x^M S_1(x) -\varphi(x,M).$$
We have
$$\gamma={\hat T_1(x)\over T_1(x)} = x^M$$
and
$$\rho= {\hat S_1(x)\over S_1(x)} = x^M-{\varphi(x,M))\over S_1(x)}.$$  

It is clear that for fixed $x$, $\varphi(x,M)$ is strictly increasing in $M$, and thus there is a unique $M^*=M^*(x)$ such that $\varphi(x,M^*(x))=0$.  Now let $\delta=\delta(\cL^c,\cL^o)$.  If $M>M^*(\delta)$, then $\varphi(\delta,M^*)>0$, which implies that $\rho(\delta)<\gamma(\delta)$.  This implies that $\hat\chi(\delta)<\gamma(\delta)\chi(\delta)=0$.  This implies that $\hat\delta>\delta$.

By inspection, $M^*(x)$ is monotone increasing in $x$, for $1\le x\le d$.  Thus for $x>\delta$ we have $M^*(x)>M^*(\delta)$.  This implies that $0=\varphi(x,M^*(x))\ge\varphi(x,M^*(\delta))$.  This means that if $M<M^*$ and $x>\delta$, then $\varphi(x,M)\le0$, which in turn implies that $\rho\ge\gamma$ and so $\hat\chi(x)\ge\gamma(x)\chi(x)>0$.  From this we conclude that $\hat\delta\le\delta$. \ \ QED

\proclaim Theorem 5.5.  Let $\cL^c,\cL^o$ be an orbit list structure, and let $\hat\cL^c,\hat\cL^o$ be the orbit list structure obtained by removing an orbit list $\cL$ from $\cL^o$ and adding $\cL$ to $\cL^c$, i.e., we move an orbit list from $\cL^o$ to $\cL^c$.  Then $\delta(\hat\cL^c,\hat\cL^o)\le\delta(\cL^c,\cL^o)$, and the inequality is strict unless $\delta(\cL^c,\cL^o)=1$.

\noindent{\it Proof.}  Without loss of generality we may assume that $\cL=\cL_1$ is
the first list in $\cL^o$, and we move $\cL$ to the first list of $\cL^c$.  Let
$T^o_1$, $T^c_1$, $S^o_1$, $S^c_1$ be defined as in (A.1).  Let $\chi$ (resp.
$\hat\chi$) be the characteristic polynomial corresponding to $\cL^c,\cL^o$ (resp.
$\hat\cL^c,\hat\cL^o$).  Then we have
$$\chi(x) = (x-d) T^o_1{\prod}'T +(x-1) S^o_1{\prod}'T + (x-1)T^o_1{\sum}'
S{\prod}''T.$$
Since the only difference between $\chi$ and $\hat\chi$ arises from the change of
$\cL$, we have
$$\gamma = {T^c_1(x)\over T^o_1(x)} = {x^{|\cL|}-1\over x^{|\cL|}}<1$$
and
$$\rho = {S^c_1(x)\over S^o_1(x)} = {S^o_1(x) +{\rm positive}\over S^o_1(x)} >1$$
for $x>0$.  It follows that $\hat\chi(x)>\chi(x)$ for $1<x<d$ and thus
$\hat\delta\le\delta$.  \ \ QED

\proclaim Theorem 5.6.  Let $\cL^c,\cL^o$ be an orbit list structure.  Suppose that $\cL=\{1,\dots,1\}$ is one of the lists inside $\cL^c$.  Let $\hat\cL^c,\hat\cL^o$ be the orbit list structure obtained by setting $\hat\cL^o=\cL^o$ and replacing the list $\{1,\dots,1\}$ in $\cL^c$ by $n$ closed lists $\{1\},\dots,\{1\}$.  Then $\delta(\hat\cL^c,\hat\cL^o)\ge\delta(\cL^c,\cL^o)$.  If $n\ge4$, then the inequality is strict.

\noindent{\it Proof.}  To fix notation, let us write $\cL=\cL_1=\{1,\dots,1\}$ be the first list of $\cL^c$, and we change $\cL$ to the first $n$ lists $\hat\cL_1=\{1\}$, \dots, $\hat\cL_n=\{1\}$ in $\hat\cL^c$.  Let $\chi$ and $\hat\chi$ be the corresponding characteristic polynomials.  Then
$$\chi(x)=(x-d)T_1{\prod}'T +(x-1)S_1{\prod}'T +(x-1)T_1{\sum}'S({\prod}''T)$$
and
$$\eqalign{ \hat\chi(x)=&(x-d)(\prod_{i=1}^n\hat T_i){\prod}^{(n)}T +(x-1)\sum_{i=1}^n S_i(\prod_{j\ne i}\hat T_j){\prod}^{(n)}T +\cr
&+ (x-1)(\prod_{i=1}^n\hat T_i){\sum}^{(n)}S({\prod}^{(n+1)}T).}$$
Thus we have
$$\gamma={\prod^n \hat T\over T_1} = {(x-1)^n\over (x^n-1)} = {(x-1)^{n-1}\over x^{n-1}+\cdots+1}$$
and
$$\rho={\sum S_i(\prod \hat T)\over S_1} = {n(x-1)^{n-1}\over n+\sum_{j=1}^{n-1} x^j}  = { (x-1)^{n-1}\over 1+\sum_{j=1}^{n-1}{n \choose j}{1\over n} x^j}.$$
Hence if $1\le x\le d$ we have $\rho\le\gamma$ and therefore $\hat\chi(x)\le\chi(x)$.  Thus $\hat\delta\ge\delta$. \ \ QED

\proclaim Theorem 5.7.  Let $f$ be an elementary mapping of ${\bf P}^d$ with $k$ singular orbits.  Then $\delta(f)\ge\delta(\cL^c,\cL^o)$, where $\cL^c=\{\{1,\dots,1\}\}$, $|\cL^c|=k$, and $\cL^o=\emptyset$.

\noindent{\it Proof.}  Let $\cL_f^c,\cL^o_f$ be the orbit list structure of $f$.  By Theorem 5.1, $\delta$ will decrease if we make all the orbits have length equal to 1.  By Theorem 5.5, $\delta$ will be also be decreased if we change all open orbit lists to closed orbit lists.  Finally, by Theorem 5.6, $\delta$ will  be decreased if we join all the orbit lists to one orbit list $\{1,\dots,1\}$.  \ \ QED

\bigskip\vfill\eject
\centerline{\bf\S 6.  Permutation mappings: I}

\noindent  Next we define the family of permutation maps.  In this section we will direct our attention to the case of the identity permutation: this is the family of mappings introduced in [BHM].  We will see that our discussion of elementary mappings applies in this case; in particular, these mappings have orbit list structure given by $\cL^c=\{\{N_1\},\dots,\{N_\ell\}\}$ and $\cL^o=\emptyset$.  Using this, we will give proofs of some conjectures from [BHM].   

Let us define
$$\cD_{j}=\{[x_{0}:\dots:x_{d}]\in\P^{d}: x_{i_{1}}=x_{i_{2}}{\rm \ for\ all\ }i_{1},i_{2}\ne j\}.$$
By $\eta_{j}(c)=(x_{0},\dots,x_{d})\in\cD_{j}$ we denote the point such that  $x_{j}=c-1$, and $x_{i}=c$ for all indices $i\ne j$.  With this notation we have $\eta_{j}(0)=e_{j}$. 
Let $a_{0},\dots, a_{d}\in\C$ be constants satisfying
$$a_{0}+a_{1}+\cdots+a_{d}=2.\eqno(6.1)$$
Let $p$ be a permutation of the set $\{0,1,\dots,d\}$, and let $P=(P_{i,j})_{0\le i,j\le d}$ be the associated permutation matrix, i.e., $P_{i,j}=\delta_{i,p(j)}$.  Let us define the $(d+1)\times (d+1)$ matrix
$$L=[\eta_{p(0)}(a_{0}),\dots,\eta_{p(d)}(a_{d})] = \pmatrix{ a_{0}&a_{1}&\dots&a_{d}\cr
\vdots&\vdots& &\vdots\cr 
a_{0}&a_{1}&\dots&a_{d}\cr}-P.$$
We set $f=L\circ J$ and define 
$$\alpha_{j}:=\eta_{p(j)}(a_{j}),\ \ \beta_{j}:=\eta_{p(j)}(1-a_{j}),\ \ \sigma_{j}=\eta_{j}(1).\eqno(6.2)$$  
It follows that $f(\Sigma_{i}\cap\cD_{i})=\sigma_{i}$.  If (6.1) and (6.2) hold, we will refer to $f=L\circ J$ as a {\it permutation mapping}.   In this case we see by the following Lemma that $f$ permutes the diagonals $\cD_{j}$ according to the permutation $p$.
\proclaim Lemma 6.1.  If (6.1) and (6.2) hold, then for each $j$ and each $c\in\C$, we have $f(\cD_{j}-e_{j})\subset\cD_{p(j)}$, and in fact: 
$$h(\eta_{j}(c))=\eta_{p(j)}(c+a_{j}-1).\eqno(6.3)$$ 

We are especially interested in the orbits $\cO(\alpha_{j})$ and $\cO(\beta_{j})$.  Since $\alpha_{j}$ and $\beta_{j}$ both belong to $\cD_{j}$, we see that these orbits can be singular (i.e., they can enter $\cI\cup\cE$) only if they end in $e_{j}$ or $\sigma_{j}$.  By the Lemma, $\cO(\alpha_{j})$ is singular exactly when one of two things happens: either
$$ a_{j}+a_{p(j)}+\dots+a_{p^{N-1}(j)} -(N-1) = 0 ,\eqno(6.4)$$
in which case the orbit of $\alpha_{j}=\eta_{p(j)}(a_{j})$ ends in $f^{N-1}\alpha_{j}=e_{p^{N}(j)}$, or
$$ a_{j}+a_{p(j)}+\dots+a_{p^{N-1}(j)} -(N-1) = 1,\eqno(6.5)$$
in which case it ends in $\sigma_{p^{N}(j)}$.
Similarly, the orbit $\cO(\beta_{j})$ is singular exactly when either 
$$ 1-a_{j}+a_{p(j)}+\dots+a_{p^{N-1}(j)} -(N-1) = 0 ,\eqno(6.6)$$
or
$$ 1-a_{j}+a_{p(j)}+\dots+a_{p^{N-1}(j)} -(N-1) = 1 .\eqno(6.7)$$

For the rest of this section we suppose that $p$ is the identity.  In this case, the only way we can have a singular orbit is in case (6.4), which becomes
$$a_{j}={N-1\over N}.\eqno(6.8)$$
It follows that all singular orbit lists consists of single orbits and are closed.  Thus the orbit list structure of our map is $\cL^{c}=\{\{N_{0}\},\{N_{1}\},\dots\{N_{k}\}\}$, $\cL^{o}=\emptyset$.  By Theorem A.1, the characteristic polynomial is given as
$$ \chi(x)=
(x-d)\prod_{j=0}^k(x^{N_j}-1)+(x-1)\sum_{j=0}^k\prod_{i\ne
j}(x^{N_i}-1). \eqno(6.9)
$$
We note from the formula for $\chi(x)$ that $f^{*}$ has a zero eigenvalue
exactly when
$\chi(0)=k+1-d=0$.  The case $k=d-1$ may be considered to be as close
as possible to the integrable case and is of particular interest.  It was
conjectured in [BHM] that the denominator of the generating function
in this case should be
$$q(x) = 1-(k+1)x-\sum_{i=1}^{k+1}(-1)^i(kx+i-1)S_i(N_1,\dots,N_n),$$
where
$S_i(N_1,\dots,N_{k+1})=\sum_{i_1<\dots<i_k}x^{N_{i_1}+\dots+N_{i_k}}$.  This is a consequence of formula (2.4) applied to the
characteristic polynomial $\chi(x)$ as given in (6.9).

Without loss of generality, we will assume that
$N_0\le N_1\le\dots\le N_k$.  If some of the $N_j$ are equal to 1, we
define $1\le\ell\le k+1$ by the condition that $N_0=\dots=N_{\ell-1}=1$ and
$N_\ell>1$.   In this case we have
$$\chi(x)=(x-1)^\ell\left[(x-\bar d)\prod_{j=\ell}^k(x^{N_j}-1)
+(x-1)\sum_{j=\ell}^k\prod_{i=\ell,i\ne j}^{k}(x^{N_i}-1)\right],\eqno(6.10)$$
where we set $\bar d=d-\ell$.
There are cases which turn out to be particularly simple:
$$\eqalign{(a)\ \  &k=d-2{\rm\ and\ }N_0=\dots=N_k=1,{\rm\ or\
}\cr
(b)\ \ &k=d-1,{\rm\ and\ }N_0=\dots=N_{k-1}=1,{\rm\ or\ }\cr
(c)\ \ &k=d-1,\ N_0=\dots=N_{k-2}=1,{\rm\ and\
}N_{k-1}=N_k=2.}\eqno(6.11)$$
In connection with conditions (6.10) and (6.11), we note that the eigenspace of $M$
corresponding to eigenvalue 1 is given by
$$\{(\tau;\alpha_0,\dots,\alpha_0;\dots;\alpha_k,\dots,\alpha_k):
(d-1)\tau+\alpha_0+\dots+\alpha_k=0\}.$$
The codimension of this space is $1+\sum_{j=0}^k(N_j-1)$.
\proclaim Theorem 6.2.  If (6.11) holds, then $d_n$ grows at most linearly.

\noindent{\it Proof.} Let $M$ be the matrix representing $f^*_X$.
In case (6.11a), $M$ is a $d\times d$ version of the matrix in (3.1), and thus
$d_n=(d-1)n+1$.

In case (6.11b) we have $\ell=k$ and $\chi=(x-1)^\ell x^{N_k}$.  The matrix
$M$ is expanded from the previous case; the lower right hand 0 is replaced
by the $N_k\times N_k$ block $\left[\matrix{0&&\cr 1&\ddots&\cr
&1&0}\right]$.  Thus $M$ has size $(1+\ell+N_k)\times (1+\ell+N_k)$ and
rank $\ell+N_k$.  The null space of $M$ is one-dimensional, and there
is a $N_k\times N_k$ block $\left[\matrix{0&&\cr 1&\ddots&\cr
&1&0}\right]$ in its Jordan canonical form.  The matrix $M-I$ is seen to
have rank $N_k$, so the rest of the Jordan canonical form consists of an
identity matrix.  It follows that $M^n=M^{N_k}$ for $n\ge N_k$.

In case (6.11c) we have $\ell=k-1=d-2$, $\bar d=2$, and
$\chi(x)=x(x+1)(x-1)^{d+1}$.  The rank of $M-I$ is seen to be 3.  This means
that the space of eigenvectors with eigenvalue 1 has dimension $d$.  Thus
the Jordan canonical from has a diagonal portion and a $2\times 2$ block
with eigenvalue 1.  The diagonal portion consists of $d-1$ ones, a zero
and a minus one.  Thus $d_{2n}$ and $d_{2n+1}$ are each a linear function of
$n$.  This completes the proof. \medskip

\proclaim Theorem 6.3.  If (6.8) holds for at least one $j$, and if (6.11) does not hold, then  $1<\delta<d$ and $\bar
d-1\le\delta\le\bar d$. 

\noindent{\it Proof.}  We saw at the end of \S4 that if (6.8) holds for at least one $j$, then $\delta<d$.  Now we show that $1<\delta$.  By Theorem 2.1, $\delta$ is the largest real zero of $\chi(x)$.   Thus we will show that if (6.11) does not hold, then $\chi$ has a zero in the interval $(1,d)$.  Let us expand $\chi$ in a Taylor series about the point
$x=1$.  We find $\chi(x)=C(x-1)^{k+1}+O((x-1)^{k+2})$, with
$$C=(1-d)N_0\cdots N_k+\sum_{j=0}^k\prod_{i\ne j}N_i.$$
We will show that if (6.11) does not hold, then $\chi<0$ on some interval
$(1,1+\epsilon)$ and thus $\chi$ will have a zero in $(1+\epsilon,d)$.

First suppose that $k\le d-2$.  Then $1-d\le-(k+1)$ and $\prod_{j\ne
i}N_i\le N_1\cdots N_k$, so
$$C\le -(k=1)N_0\cdots N_k+(k+1)N_1\cdots N_k.$$
This gives $C<0$ unless $N_0\cdots N_k=\prod_{i\ne j}N_i$, in which case
$N_0=N_1=\dots=N_k=1$.  This is case (6.11a).

Now suppose $k=d-1$, let $\ell$ be as above, and factor $\chi=(x-1)^\ell
q(x)$ as in (6.10).  Expanding $q(x)$ about $x=1$, we obtain $q(x)=\tilde
C(x-1)^{k-\ell+1}+O((x-1)^{k-\ell+2})$,
where
$$\tilde C=(1-\bar d)N_\ell\cdots N_k+\sum_{j=\ell}^k{\prod_{\ell\le
i\le k}}'N_i,$$
where ${\prod}'$ means that the product is taken over $i\ne j$.  Since
$N_\ell\ge2$, we have
$$\tilde C\le 2(\ell+1-d)N_{\ell+1}\cdots N_k+(kl-\ell+1)N_{\ell+1}\cdots
N_k =(\ell+2-d)N_{\ell+1}\cdots N_k.$$
Now we may assume $\ell\le k-1$, for otherwise we are in case (6.11b).  Thus
$\tilde C<0$ unless $\ell=d-2$, and thus $d-1=k$.  We have already
handled the case $\ell=k=d-1$.  Thus we have $\tilde C<0$ unless
$\ell=k-1=d-2$.  By our formula, then,
$$\tilde C=(d-2+1-d)N_{k-1}N_k+N_{k-1}+N_k=-N_{k-2}N_k+N_{k-1}+N_k,$$
which is strictly negative unless $N_{k-1}=n_k=2$, which is case (6.11c). 
This completes the proof.\medskip

\proclaim Theorem 6.4.  If (6.11) does not hold, then $\delta$ is a simple
eigenvalue for
$f^*$.  If, in addition, $\bar d\ge3$, then $\delta$ is the unique root
of $\chi$ in the interval $[2,d]$.

\noindent{\it Proof.}  By (6.10) we may assume that $\bar d=d$, which is to say
that $N_j\ge2$ for all $j$.  We will suppose that $x>1$ is a zero of
$\chi$, and we will show that for such a zero we have $\chi'(x)>0$.  If we
divide $\chi$ by $(x^{N_j}-1)$, the condition that $\chi(x)=0$ is
equivalent to $$(x-d)\prod_{i\ne j}(x^{N_i}-1)= -{x-1\over
x^{N_j}-1}\prod_{i\ne j}(x^{N_i}-1)-(x-1)\sum_{i\ne
j}{\prod}''(x^{N_\ell}-1),$$
where ${\prod}''$ indicates a product over all $\ell$ distinct from $i$ and
$j$.  In order to compute $\chi'(x)$ we first use the product rule and then
we substitute the identity above to obtain
$$\chi'(x)=\left(\prod_j(x^{N_j}-1)\right)\left[1+\sum_j\left\{
1-(x-1){N_jx^{N_j-1}\over x^{N_j}-1}\right\}{1\over x^{N_j}-1}\right].$$
We abbreviate this as
$$\chi'(x)=\left(\prod_j(x^{N_j}-1)\right)\left[1+\sum_{j}\varphi_j(x)\right],$$
where each $\varphi_j$ has the form
$$\varphi(x)={1\over x^N-1} +(1-x){Nx^{N-1}\over (x^N-1)^2}.$$
Since the product in the formula for $\chi'$  is strictly positive for
$x>1$, it suffices to show that $\varphi_j(x)>0$ for each $j$.  In fact, we
have $\varphi(x)>0$ for all $x\ge2$ and $N\ge2$.
For this, we note that $\lim_{x\to\infty}\varphi(x)=0$, and we show that
$$\varphi'(x)={Nx^{N-2}\over (x^N-1)^3}\left[ 3(N-1)x^{N+1}
-(3N-1)x^N-(N-3)x+(N-1)\right]<0.$$
This is equivalent to showing that the expression in square brackets is
positive for all $N\ge2$ and $x\ge 2$.  This is elementary, and so we
conclude that
$\chi'(x)>0$ for every zero of $\chi$ in the interval $[2,d]$.  Thus there
can be no more than one zero in $[2,d]$.

If $d\ge3$, then as was observed above, $2\le d-1\le\delta\le d$.  Thus
$\delta $ is a simple zero of $\chi$.  If $d=2$, then $1\le\delta\le 2$,
and so the arguments above do not apply directly.  However, the case $d=2$
may be broken into three subcases (1) $k=0$, $N_0\ge2$, (2) $k=1$,
$N_0=2<N_1$, and (3) $k=1$, $3\le N_0\le N_1$.  The computations are
similar to what we have done already, so we omit the details.

\bigskip\centerline{\bf \S7. Examples}

\medskip\noindent A number of mappings of the form $L\circ J$ have arisen in the mathematical physics literature.  Let us show how the preceding discussion may be applied to yield the degree complexity of these maps.  The third example will lead us to some
non-elementary maps, and our treatment of them will foreshadow the technique we
use in \S8.

\noindent{\bf Example 7.1.}  We consider the (families of) matrices:
$$A_{1} = \pmatrix{1&2 &1\cr 1&0&-1\cr 1&-2&1\cr}, \ A_{2}= \pmatrix{ 2  &-1 + q^2&-1 + q^2\cr 2 &-1 + q  &  -1 - q\cr  2  &  -1 - q  &  -1 + q\cr} $$
$$A_{3} = \pmatrix { 1 -{\left( 1 + \ell  \right) \,\left( -1 + 2\,\ell  \right) \over \ell \,\left( 1 + 2\,\ell  \right) }  & -{1\over \ell }  &
   {3 + 2\,\ell \over 1 + 2\,\ell } \cr
 {1 - 2\,\ell   \over  \ell \,\left( 1 + 2\,\ell  \right) }&
   1 - {1\over  \ell \,\left( 1 + \ell  \right) }&
   {3 + 2\,\ell  \over \left( 1 + \ell  \right) \,\left( 1 + 2\,\ell  \right) }\cr 
{-1 + 2\,\ell \over  1 + 2\,\ell }&
    {1  \over  1 + \ell }&
   1 -  {\ell \,\left( 3 + 2\,\ell  \right) \over  \left( 1 + \ell  \right) \,\left( 1 + 2\,\ell  \right) }\cr}$$
  We set $f=A_{j}\circ J$ for $j=1,2,3$.  (The case of matrix $A_{1}$ arises, for instance,  in [BMV], $A_{2}$ in [BMV] and  [V], and $A_{3}$ is found in [R3]; and also in [R1] for the special case $\ell=1$.)  In each case we have
$$\Sigma_{j}\to\alpha_{j}\to e_{j},\ \ \ 0\le j\le 2.$$
Thus the orbits $\cO_{j}$ are singular, and $|\cO_{j}|=2$ for $0\le j\le 2$.  
In other words, the singular orbit list structure is
$\cL^{c}=\{\{2\},\{2\},\{2\}\}$ and $\cL^{o}=\emptyset$.  If we write
$f^{*}_{X}$ according to equations (4.1) (4.2) and (4.3), we find that
$(f^{*}_{X})^{6}$ is the identity.  Thus $\d(f^{n})$ is bounded since it  is
a periodic sequence of period 6.

\medskip\noindent {\bf Example 7.2.}  We consider the matrices:
$$B_1={1\over 9}\pmatrix {1&-8&16\cr
-2&7&4\cr
4&4&1},\ \ B_2={1\over 10}\pmatrix{1&22&77\cr
1&-8&7\cr
1&2&-3\cr}.$$
(We have taken $B_{1}$ from [R2] and $B_{2}$ from [V].)  Let $g_j=B_j\circ J$.  It follows that for both $g_1$ and $g_2$ we have:
$$\Sigma^*_0\to *\to e_0.$$
In the case of $g_1$, the point $[1:1:1]$ is a parabolic fixed point, and the
orbits $\cO_{1}$ and $\cO_{2}$ are in the attracting basin of $[1:1:1]$, so they are both nonsingular.  Similarly, in the case of
$g_2$, the orbits $\cO_{1}$ and $\cO_{2}$ are in the basin of an atracting
2-cycle and thus are both nonsingular.  We conclude that the
singular orbit list structure for both $g_1$ and $g_2$ is $\cL^c=\{\{2\}\}$,
$\cL^o=\emptyset$, and so
$\delta(g_1)=\delta(g_2)=(1+\sqrt5)/2$ is the largest root of $x^2-x-1=0$.

\medskip\noindent{\bf Example 7.3.}  Consider the family of matrices:
$$C(q) =\pmatrix{ 2 & -1 - q^2 & -1 - q^2 \cr 2 & -1 + q & -1 - q\cr
 2 & -1 - q  & -1 + q \cr}$$
with $q\ne0$, which is considered in [V].  If we set $h=C(q)\circ J$, then $\Sigma_{0}\to\alpha_{0}\to e_{0}$ for all  $q$.  The natures of the orbits $\cO_{1}$ and $\cO_{2}$, however, are dependent on $q$.  For generic $q$, the orbits $\cO_{1}$ and $\cO_{2}$ are nonsingular, so we have $\cL^{c}=\{\{2\}\}$, $\cL^{o}=\emptyset$, so $\delta=(1+\sqrt5)/2$ as in Example 2.

Now let us show what happens in the singular cases.  Our purpose here is to show how the methods of \S4 can be used to treat the different cases that can arise.  First let us handle the most singular cases:
\medskip\noindent{\bf Case $q=-1$.} $\alpha_{0}\to e_{0}$, $\alpha_{1}=e_{1}$, $\alpha_{2}=e_{2}$.  In this case we have $\cL^{c}=\{\{2\},\{1\},\{1\}\}$, $\cL^o=\emptyset$, and the degrees are periodic of period 3.

\medskip\noindent{\bf Case $q=1$.} $\alpha_{0}\to e_{0}$, $\alpha_{1}=e_{2}$,
$\alpha_{2}=e_{1}$.  In this case we have $\cL^{c}=\{\{2\},\{1,1\}\}$, $\cL^o=\emptyset$, and  the degrees are periodic with period 6. 

In every case, we have $\alpha_{0}\to e_{0}$, so we pass to the map $h:X\to X$ obtained by blowing up the orbit $\cO_{0}=\{\alpha_{0},e_{0}\}$.  Now observe:
$${\rm\   If\ \ } h[a:b:c]=[a':b':c'],{\rm\  then\ \ } h[a:c:b]=[a':c':b'].\eqno(7.1)$$  
It follows from (7.1) that $h^{n}\alpha_{1}\in\cI\cup\cE$ if and only if $h^{n}\alpha_{2}\in\cI\cup\cE$.  It follows that we may proceed by induction on $n:=|\cO_{1}|=|\cO_{2}|$.  Let us define
$$S':=\{x_{0}-x_{2}=0\}, \ \ \ S'':=\{(1+q)x_{1}+(-1+q)x_{2}=0\}.$$
We have $h: S'\leftrightarrow S''$.  Since $\alpha_{1}\in S''$, it follows that $h^{n}\alpha_{1}\in S'$ when $n$ is odd and $h^{n}\alpha_{1}\in S''$ when $n$ is even.   By (7.1), we have 
$$\{x_{0}-x_{1}=0\}\leftrightarrow \{(1+q)x_{2}+(-1+q)x_{1}=0\},$$
so an analogous discussion applies to the orbit of $\alpha_{2}$.

\medskip\noindent{\bf Case $n=|\cO_{1}|=|\cO_{2}|$ is even.}  If $n$ is even, then
$h^{n-1}\alpha_{j}\in\cE\cup\cI$.  If $h^{n-1}\alpha_{j}\in \cI$, then $h^{n-1}\alpha_{j}\in S'\cap\cI$,
and we have $h^{n-1}\alpha_{j}=e_{j}$ for $j=1,2$.  Thus our orbit structure is
$\cL^{c}=\{\{2\},\{n\},\{n\}\}$ and $\cL^{o}=\emptyset$. If $n=2$, then $d_{n}$ is periodic of period 6, as in Example 7.1.  If $n\ge6$, then the degree complexity $\delta_n$ is the
largest root of the polynomial $x^{n+2}-x^{n+1}-x^n+x^2+x-1$.  

The other possibility is that $h^{n-1}\alpha_{1}=[1:0:1]\in S'\cap\Sigma_{j}$.  In this case, we have $h^{n}\alpha_{1}=\alpha_{1}$.  By (7.1), a similar argument applies to $\alpha_{2}$. Thus $\alpha_{1}$ and $\alpha_{2}$ are periodic, and the orbits $\cO_{1}$ and $\cO_{2}$ are essentially nonsingular.  Thus we have $\cL^{c}=\{\{2\}\}$, and we have $\delta=(1+\sqrt5)/2$ as in Example 7.2.
\medskip\noindent{\bf Case $n=|\cO_{1}|=|\cO_{2}|$ is odd.}  In this case we have $h^{n-1}\alpha_{1}\in S''\cap\cE\cup\cI$.  We cannot have $h^{n-1}\alpha_{1}\in\cI\cap S''$, since the point $e_{0}$ has been blown up, and we can only reach the blow-up fiber $\cF_{e_{0}}$ through the fiber $\cF_{\alpha_{0}}$, and we can reach $\cF_{\alpha_{0}}$ only through $\Sigma_{0}$.  Thus we must have $h^{n-1}\alpha_{1}=[0:1-q:1+q]\in S''\cap\Sigma_{0}$.
Let us consider this point as the endpoint of the curve $t\mapsto[t:1-q:1+q]$ as $t\to0$.  Thus $h$ maps
this to the curve 
$$t\to (1,1,1) +t({1+q^{2}},1-q^{2}, {1+q^{2}})+O(t^{3})$$ which lands at a point of the fiber $\cF_{\alpha_{0}}$, and then to the curve 
$$t\to(1,0,0) +t(0,{ 1-q},{1+q})+t^{2}(1+q^{2},2(1+q)^{-1},2 (1-q)^{-1}) +O(t^{3}),$$
 which lands at a point of the fiber $\cF_{e_{0}}$.  The next image of this curve lands at 
$$\beta_{1}:=h^{n+2}\alpha_{1}= [1+q^{2}:1-q^{2}:1+q^{2}] \in{\bf P}^{2}-\{\alpha_{0},e_{0}\},$$
and we are back to a ``normal'' point of $S'$.

At this stage there are two possibilities.  First, it is possible that $\cO(\beta_{1})$ is nonsingular. We conclude, then that $h:X\to X$ is 1-regular, and we have $\delta=(1+\sqrt5)/2$ as in Example 2.
The other possibility is that $\cO(\beta_{1})$ is singular.  This means that $h^{j}\beta_{1}\in X-(\cI\cup\cE)$ for $0\le j< j_{1}$, and $h^{j_{1}}\beta_{1}\in\cI\cup\cE$. First, we see that $h^{j_{1}}\beta_{1}$ cannot be in $S''\cap\cE$.  For in this case we must have $h^{j_{1}}\beta_{1}=[0:1-q:1+q]$ as before.  But this is not possible since we have remained inside points where $h$ is a diffeomorphism.  On the other hand, if we have $h^{j_{1}}\beta_{1}\in S'$, then we must have $h^{j_{1}}\beta_{1}=[1:0:1]$.  We have $h^{j_{1}}\beta_{1}\in S'$ if $j_{1}$ is even.   Thus $h^{j_{1}+1}\beta_{1}=\alpha_{1}$, and $\alpha_{1}$ is periodic.  A similar argument shows that both $\alpha_{1}$ and $\alpha_{2}$ are periodic in this case.  Thus $\cO_{1}$ and $\cO_{2}$ are both essentially nonsingular, and we are in the case of Example 2 again.

\noindent{\bf Sub-case $\cO(\beta_{1})$ is singular.}   The other possibility is that $\cO(\beta_{1})$ ends at the point $e_{1}\in S'$.  In this case, $\cO(\beta_{2})$ also ends at $e_{2}$, and $|\cO(\beta_{1})|=|\cO(\beta_{2})|$.  This sub-case is not elementary, and here we must perform a second series of blow-ups.   Let $\tilde\cO(\alpha_{1})$ denote the orbit in $X$, starting with $\alpha_{1}$.  Figure 1 shows $\tilde\cO(\alpha_{1})$ and  $\tilde\cO(\alpha_{2})$ in the space $X$.   On the top row, the portion $\alpha_{1}\to(*)_{1}\to [0:1-q:1+q]$ is the orbit $\cO_{1}=\cO(\alpha_{1})$, and $(*)_{1}$ indicates the points in the middle of the orbit.  The image $h[0: 1-q: 1+q]$ is indicated by the subscript $\alpha_{0}$ (base point) and fiber coordinate $[1+q^{2}:1-q^{2}:1+q^{2}]\in\cF(\alpha_{0})$.  We will use the notation
$$\tau_{0,1}=\left[\matrix{0\cr 1-q\cr 1+q}\right]_{ e_{0}}, \ \ \tau_{0,2} = \left[\matrix{0\cr 1+q\cr 1-q}\right]_{e_{0}}$$
for the points of the orbit that are in the fiber $\cF(e_{0})$.  The bottom row shows the orbit of $\Sigma_{0}$, which contains the two points of exceptional fibers $\cF(\alpha_{0})$ and $\cF(e_{0})$.

$$\matrix{ \Sigma_{1}  \to  \alpha_{1}  \to  (*)_{1}\to & \left[\matrix{0\cr 1-q\cr 1+q}\right]  \to &\left[\matrix{1+q^{2}\cr 1-q^{2}\cr 1+q^{2}}\right]_{ \alpha_{0}} \to &\tau_{0,1} & \to  \beta_{1}  \to  (**)_{1}  \to  e_{1} \cr
\Sigma_{2}  \to  \alpha_{2}  \to  (*)_{2} \to & \left[\matrix{0\cr 1+q\cr 1-q}\right] \to & \left[\matrix{1+q^{2}\cr 1+q^{2}\cr 1-q^{2}}\right]_{\alpha_{0}}  \to & \tau_{0,2}&  \to  \beta_{2}  \to  (**)_{2}  \to  e_{2} \cr
&&&\phantom{0}&\cr
& \Sigma_{0} \ \ \ \ \  \to & \cF(\alpha_{0})\ \ \ \ \ \to & \cF(e_{0}) & \cr }$$
\centerline{Figure 1}
\smallskip
 We let $\pi:X_{2}\to X$ denote the space obtained from $X$ by blowing up the points of $\tilde\cO(\alpha_{1})\cup\tilde\cO(\alpha_{2})$.  Let $h_{2}:X_{2}\to X_{2}$ denote the induced birational map.  We see that in $X_{2}$ the curves $\Sigma_{1}$ and $\Sigma_{2}$ are no longer exceptional, so $h_{2}$ is 1-regular.  Let us determine the mapping on cohomology.  As a basis for $H^{1,1}(X_{2})$ we take $H_{X_{2}}$, together with all the fibers indicated in Figure 1.  For instance, we take the fibers $\cF(p)$ for all the points $p\in (**)_{1}$.  We see that under $h_{2}^{*}$ we have:
 $$\cF(e_{0})\to\cF(\alpha_{0})\to\{\Sigma_{0}\}_{X_{2}}\eqno(7.2)$$
 and
 $$\cF(e_{1})\to \cF(p)\to\cdots\to\cF(p) \to \cF(\beta_{1})\to\cdots\to\cF(\alpha_{1})\to\{\Sigma_{1}\}_{X_{2}},\eqno(7.3)$$
where the $p$ are points in the $(**)_{1}$ portion of the orbit, taken in inverse order.  Similarly, we note that the centers of blow-up inside $\Sigma_{0}$ are $e_{1}$, $e_{2}$, $[0:1-q:1+q]$, and $[0:1+q,1-q]$.  Thus we have
 $$\{\Sigma_{0}\}_{X_{2}}=H_{X_{2}} - \cF(e_{1})-\cF(e_{2}) - \cF\left[\matrix{0\cr 1-q\cr 1+q}\right] -\cF\left[\matrix{0\cr 1+q\cr 1-q}\right].\eqno(7.4)$$
 On the other hand, $e_{0}$, $e_{2}$ are the base points in $\Sigma_{1}$.  The closure of $\Sigma_1$ in $X_2$ intersects the fiber $\cF(e_0)$ at the point with (fiber) coordinate $[0:0:1]$. Since $q^2 \ne 1$, it follows that this is distinct from the base points $\tau_{0,1}$ and $\tau_{0,2}$ of the second level blow up. Thus we have
 $$\{\Sigma_{1}\}_{X_{2}} = H_{X_{2}} - \cF(e_{0}) - \cF(e_{2})-\cF(\tau_{0,1})-\cF(\tau_{0,2})\eqno(7.5)$$
 and
 $$\{\Sigma_{2}\}_{X_{2}} = H_{X_{2}} - \cF(e_{0}) - \cF(e_{1})-\cF(\tau_{0,1})-\cF(\tau_{0,2}).\eqno(7.6)$$
 Finally, we must evaluate $h^{*}_{2}H_{X_{2}}$.  If $\cH$ is a general hypersurface in $\P^{2}$, then $\{e_{0},e_{1},e_{2}\}\subset h^{-1}\cH=J^{-1}L^{-1}\cH$.  We have seen that $J_{X}$ is nonconstant on the fibers $\cF(e_{j})$, so $e_{0},[1+q^{2}:1-q^{2}:1+q^{2}]$ will not be contained in $J_{X}^{-1}\cH$ for generic $\cH$.  We conclude that
 $$h^{*}_{2}H_{X_{2}} = 2H_{X_{2}}-\cF(e_{0})-\cF(e_{1})-\cF(e_{2})-\cF(\tau_{0,1})-\cF(\tau_{0,2}).\eqno(7.7)$$
Equations (7.2--7) serve to define the linear transformation $f^{*}$.  Assuming that $|\cO(\alpha_{1})|=|\cO(\beta_{1})|=n$, the characteristic polynomial of this transformation turns out to be the same as the characteristic polynomial corresponding to the elementary map $\cL^{c}=\{\{2\},\{2n+2\},\{2n+2\}\}$, $\cL^{o}=\emptyset$.
 
 \noindent{\bf Observed cases.}  The first few cases with $n$ even occur for $n=2$ if $3+q^2=0$, $n=6$ if  $5+10q^{2}+q^{4}=0$, $n=10$ if $7+ 35q^{2}+21q^{4}+q^{6}=0$.  In all of these cases the orbit $\cO_j$ ends with $e_j$.  The first few cases with $n$ odd occur for $n=1$ if $1+q^2=0$, $n=3$ if $1+3q^2=0$, $n=5$ if $q^{4}+6q^{2}+1=0$, $n=7$ if $1 + 10q^{2}+5q^{4}=0$, $n=9$ if $1+14q^2+q^4=0$, $n=11$ if $1+21q^{2}+35q^{4}+7q^{6}=0$.  In all of the odd cases, the orbit $\cO(\beta_j)$ is singular, and $|\cO(\alpha_j)|=|\cO(\beta_j)|$.  In both the even and odd cases,  $\delta_n$ is less than the generic case $(1+\sqrt5)/2$, and we see from the defining equations that $\delta_n\to(1+\sqrt5)/2$ as $n\to\infty$.

\bigskip\centerline{\bf\S8. Permutation mappings: Orbit collision, orbit separation.}
\medskip
\noindent  Here we continue our discussion of permutation mappings.  Like Noetherian mappings, the permutation mappings have the form $f=L\circ J$.  As before, the key to understanding these mappings is understanding what happens with the orbits of the points $\alpha_j:=f\Sigma_j^*$.  If such an orbit is singular, it ends at a point $e_k$ or a point $\sigma_k\in\Sigma_k^*$.  The first case corresponds to elementary behavior and has been treated above.  The second case, corresponding to (6.4--7), is an example of non-elementary behavior; in this case the orbit $\cO(\alpha_j)$ ``joins'' the orbit $\cO(\alpha_k)$.  We refer to this as an orbit collision.  Further orbit collisions are also possible, with $\cO(\alpha_k)$ joining $\cO(\alpha_m)$, etc.  Our interest in \S8 is to show that the method of regularization can be applied to the case of orbit collisions.  This leads us to perform multiple blow-ups over a fixed base point, which provides a new manifold in which these orbits are separated, and the induced map is 1-regular.

Let $\cS$ denote the set of all orbits $\cO(\alpha_{j})$, $\cO(\beta_{j})$, $0\le j\le d$, which are singular.  There are four possible types of singular orbits in $\cS$: $\alpha e$, $\alpha\sigma$, $\beta e$ and $\beta\sigma$, depending on the type of starting point and the type of ending point.  Now we define admissible chains of singular orbits.  The admissible chains of the first generation are the singular chains starting with an $\alpha$ and ending with an $e$.  We denote the chains of the first generation by $\cC^{1}$.   Now let us proceed inductively, assuming that we have defined the admissible chains $\cC^{j}$ at generation $j$.  An admissible chain of generation $j+1$ will be a finite sequence of singular orbits of $\cS$ which has the following form: $\cS\cC\cS\cC\cS\dots\cC\cS$, which means that we start and end with orbits of $\cS$ and in the middle, we alternate between $\cS$ and $\cC=\cC^1\cup\dots\cup\cC^{j}$. By convention $\cC^{j+1}$ is disjoint from $\cC^1\cup\dots\cup\cC^{j}$. In addition, the sequence must obey the following rules: 
The first orbit starts with an $\alpha$;  the last orbit ends with an $e$, and the permissible transitions between $\cS$ and $\cC$ are $e_{i}\to\beta_{p(i)}$ and $\sigma_{i}\to\alpha_{p(i)}$.   In other words, suppose that $\cO'$ is an orbit from $\cS$ which is followed by a chain $\cO''\dots\cO'''\in\cC$.  If $\cO'$ ends with $e_{\ell}$, then $\cO''$ must begin with $\beta_{p(\ell)}$; and if $\cO'''$ ends with $\sigma_{\kappa}$, and if $\cO'''$ is followed by an orbit $\cO''''\in\cS$, then $\cO''''$ must begin with $\alpha_{p(\kappa)}$.  The process of constructing chains is finite, so there is a maximum generation $\kappa$ that can occur.  Thus $\cC^{1}\cup\cC^{2}\cup \dots\cup\cC^{\kappa}$ is the set of admissible chains.    

We illustrate this with an example, which is sketched in Figure 2.  This corresponds to the cyclic permutation $p=(1,2,3,\dots,N)$ with $N$ greater than 14.  The singular orbits $\cS$ are $\{\alpha_{1},\sigma_{2}\}$, $\{\alpha_{3},\sigma_{4}\}$, $\{\alpha_{5},e_{6}\}$, $\{\beta_{7},e_{8}\}$, etc.  These are inside the bottom row of Figure 1.  To conserve space, we have constructed all of these orbits to have (minimal) length 2.  The chains of the $j$th generation may be read off from the $\cC^{j}$ row of the matrix by joining adjacent dots, moving from left to right.  Thus $\cC^{1}=\{\{\alpha_{5},e_{6}\}, \{\alpha_{11},e_{12}\}\}$ consists of two chains, and $\cC^{2}=\{\alpha_{3},\dots, e_{8}\}$ and $\cC^{3}=\{\alpha_{1},\dots,e_{14}\}$ each contain one chain.
\medskip
$$\matrix{ \cC^{3}& & \cdot &  \cdot & \cdot & \cdot & \cdot & \cdot & \cdot & \cdot & \cdot & \cdot & \cdot & \cdot & \cdot&\cdot \cr
\cC^{2} && &&  \cdot & \cdot & \cdot & \cdot & \cdot & \cdot & &&&&&\cr
\cC^{1} && &&&&  \cdot & \cdot & &&&&  \cdot & \cdot &&\cr
\cr
&& \alpha_{1} &\sigma_{2}& \alpha_{3} &\sigma_{4} &\alpha_{5}&e_{6}&\beta_{7} &e_{8} &\beta_{9}& \sigma_{10} &\alpha_{11} &e_{12}&\beta_{13}&e_{14}\cr} $$
\centerline{Figure 2.  Singular Chains}

In order to construct a 1-regularization of $f$, we perform multiple blow-ups of points of $\P^d$, determined by the structure of the chains.  We define the {\it height} of a point $p\in\P^d$, written $h(p)$, to be the number of chains $\gamma$ that contain $p$.   Note that if $\gamma\in\cC^{j}$ is a singular chain, the height $h(p)$ changes by at most 1 as we step forward from one point $p\in\gamma$ to the next one.  Let $\cP$ denote the set of points $p$ of $\P^d$ which occur in singular chains.  Thus $\cP=\{p:h(p)>0\}$.   Observe that $\alpha_j,\beta_j\in\cD_j$, and thus $\cP\subset\bigcup_{j=0}^d\cD_j$.   

We will define the space $X$ by blowing up $h(p)$ times over each $p\in\cP$.  Let $\pi^1:X^1\to\P^d$ denote the space obtained by blowing up $\P^d$ at each $p\in\cP$.   As we construct manifolds $\pi^k:X^k\to X^{k-1}$, it will be convenient to let $\cD_j\subset X^k$ denote the strict transform of $\cD_j$, i.e. the closure of $(\pi^k)^{-1}(\cD_j-\{p\})$ in $X^k$.   The strict transform of $\cD_j$ in $X^1$ intersects $\cF^1(p)$ transversally at a point $p^1:=\cD_j\cap\cF^1(p)$.  We now construct $\pi^2:X^2\to X^1$ by blowing up all the points $p^1\in\cF^1(p)$ for which $h(p)>1$.  Let $\cF^2(p)=(\pi^2)^{-1}(p^1)$ denote the new fiber.  To simplify our notation we write $\cF^1(p)$ for the strict transform of $\cF^1(p)$ inside $X^2$.  This abuse of notation causes no problem because $\cF^1(p)\cap\cF^2(p)$ has codimension 2. Since $p^1\in\cD_j$, it follows that $\cD_j$ intersects $\cF^2(p)$ transversally at a point $p^2$.  We continue the blow-up process at the points $p^2$ for which $h(p)>2$.  We continue in this way until we reach the maximum value of $h$; thus we construct the space $X$.

It follows that over every point $p\in\cP$, we have exceptional fibers $\cF^j(p)$, $1\le j\le h(p)$.   For simplicity of notation, we let $\cF^j(p)$ denote its corresponding class in $H^{1,1}$.  These cohomology classes, together with the class $H_X$ of a hyperplane,  generate $H^{1,1}(X)$.   We find it convenient to use the notation $\hat\cF(p)=\sum_{j=1}^{h(p)} \cF^j(p)$.

Let us describe how the induced map $f_X:X\to X$ maps the various exceptional fibers.  Let us start with a singular chain $\cO=\{\alpha_j,\dots,e_k\}$ of the first generation.  By \S3, the first and last maps in the sequence
$$\Sigma_{p^{-1}(j)}\to\cF^1(\alpha_j)\to\cdots\to\cF^1(e_k)\to L\Sigma_k\eqno(8.1)$$ 
have (maximal) generic rank $d$.
Since $f$ is locally biholomorphic at each point of $\cO-\{e_k\}$, the rest of the maps are biholomorphic in a neighborhood of the fibers.  In particular, none of the hypersurfaces in (8.1) is exceptional.  

Next, consider a singular chain of the 2nd generation.  We may suppose that the
chain has the form $\cO'\cO''\cO'''$, where
$\cO'=\{\alpha_i,\dots,\sigma_{p^{-1}(j)}\}$, $\cO'' =\{\alpha_j,\dots,e_k\}$,
and $\cO'''=\{\beta_{p(k)},\dots,e_\ell\}$.  We will now claim that $f_X$
induces full rank mappings
$$\eqalign{\Sigma_{p^{-1}(i)}\to\cF^1(\alpha_i) &
\to\cdots\to\cF^1(\sigma_{p^{-1}(j)})\buildrel
1\over\to\cF^2(\alpha_j)\buildrel 2\over\to\cdots\cr &\cdots\buildrel
2\over\to\cF^2(e_k)\buildrel
3\over\to\cF^1(\beta_{p(k)})\to\cdots\to\cF^1(e_\ell)\to
L\Sigma_\ell.}\eqno(8.2)$$ 
We need to discuss the arrows marked with ``1,'' ``2,'' and
``3.''  The arrows marked ``2'' come about because $f$ maps $\cD_j$ to
$\cD_{p(j)}$, and so $f$ maps $p^1$ to $(fp)^1$.  Since $f$ is locally
biholomorphic at $p$, it follows that $f$ is biholomorphic in a neighborhood of
$\cF^2(p)$.  

To analyze maps ``1'' and ``3,'' we give a local coordinate system on $X$.   Let us start with an affine coordinate chart about $p=0$.  For $v\in\C^d$ we  let $\gamma_v: t\mapsto tv$.  We identify $v$ with a point in the fiber  $\cF^1(p)$.   The point $p^1\in\cF^1(p)$ corresponds to the landing point of $\gamma_{\bf 1}$, where  ${\bf 1}=[1:\dots:1]$.  For $w\in\C^d$, we consider
the family of curves
$\gamma_{{\bf 1},w}:t\mapsto t{\bf 1}+t^2w$, which are tangent to
$\cD_i$ at $p$.  To consider $w$ as a coordinate of the fiber $\cF^2(p)$, we note that if we reparametrize the curve by the change of variable $t\leftarrow t+at^2$, we modify $w$ by adding a multiple of ${\bf 1}$.   Using the local coordinate $w$, we may identify a point of $\cF^2(p)$ with the quotient $\C^d/\langle{\bf 1}\rangle$.

Now for
map ``1,'' let us consider the map $J$ in a neighborhood of
$\sigma_0=[0:1:\dots:1]$.  A point of the local coordinate chart $v_0\ne0$ of 
$\cF^1(\sigma_0)$ is identified with the landing point of $\gamma_v: t\mapsto
\sigma_0+t(1,v_1,\dots,v_d)$.  Under $J$, this is mapped to
$$t\mapsto e_0+t\sigma_0 + t^2(0,-v_1,\dots,-v_d)+O(t^3).\eqno(8.3)$$ 
It follows that $J$ is injective on a dense open set of $\cF^1(\sigma_0)$, so map ``1'' has full rank.  

Since $J$ is an involution, it maps the curve (8.3) back to $\gamma_v$, so $J$ induces a full rank map
$\cF^2(e_0)\to\cF^1(\sigma_0)$.  Thus map ``3'' has full rank.   We conclude that none of the hypersurfaces in
(8.2) is exceptional.  Continuing by induction to the higher blow-up fibers, we have the following:

\proclaim Proposition 8.1.  The induced map $f:X\to X$ is 1-regular.

Let us describe the induced mapping $f^*$ on cohomology.  By the discussion above, we see that for any singular chain $\gamma$, there is a sequence of mappings through hypersurfaces $\cF^j(p)$: 
$$\Sigma_{s(\gamma)}\to \cF^1(\alpha_{p(s(\gamma))})\to\cdots\to\cF^1(e_{\omega(\gamma)}).$$
If $\gamma'$ and $\gamma''$ are distinct chains, then the endpoints are distinct $s(\gamma')\ne s(\gamma'')$, and $\omega(\gamma')\ne\omega(\gamma'')$.  The mapping on cohomology is given by:
$$\cF^1(e_{\omega(\gamma)})\to \cdots\to\cF^1(\alpha_{p(s(\gamma))})\to\Sigma_{s(\gamma)}.\eqno(8.4)$$
We see that every $\cF^j(p)$ is contained in a unique singular chain $\gamma$, so (8.4) tells how $f_X^*$ acts on each $\cF^j(p)$.  

Now let us write $\Sigma_0$ with respect to our basis.  The only centers of blow up inside $\Sigma_0$ are $\{\sigma_0,e_1,\dots,e_d\}=\cP\cap\Sigma_0$.  Now each $\cD_0$ is the line connecting $\sigma_0$ to ${\bf 1}\notin\Sigma_0$.  Thus $\cD_0$ intersects $\Sigma_0$ transversally.  Thus $\cD_0\cap\cF^1(\sigma_0)$ is disjoint from $\Sigma_0\cap\cF^1(\sigma_0)$.  Similarly, for $j\ne0$, $\cD_j$ is the line connecting $e_0$ to  ${\bf 1}\notin\Sigma_0$.  Thus $\cD_j$ intersects $\Sigma_0$ transversally.  Thus $\cD_j\cap\cF^1(e_0)$ is disjoint from $\Sigma_0\cap\cF^1(\sigma_0)$.  We conclude that
$$\{\Sigma_0\}_{X}=H_X-\sum_{j\in\cP-\{e_0\}}\hat\cF(e_j) -\hat\cF(\sigma_0)\in H^{1,1}(X).\eqno(8.5)$$
By a similar argument, we have
$$f^*H_X=d\cdot H_X +(1-d)\sum_{e_j\in\cP}\hat\cF(e_j).\eqno(8.6)$$
It follows that $f_X^*$ is given by (8.4--6).

\bigskip
\centerline{\bf  Appendix:  Characteristic Polynomial}

\medskip\noindent Let $\cL=\{\cL_1,\dots,\cL_\mu\} =\{\{N_{1,1},
\dots, N_{1,\ell_1}\},\{N_{2,1}, \dots, N_{2,\ell_2}\}, \dots
\{N_{\mu,1}, \dots, N_{\mu,\ell_{\mu}}\}\}$ denote the set of
lists of lengths of orbits inside orbit lists.  Let us fix an
orbit list $\cL_i$ and let $M$ denote a subset of indices
$\{1,\dots,\ell_i\}$.  We define 
$$|\cL_i|_M=\sum_{j\in M}N_{i,j}.$$  
If $\cL_i=\cL^o_i$ is an open orbit list, we define 
$$T^o_i(x)=x^{|\cL^o_i|},{\rm\ \ and \ \ }S^o_i(x) = \sum_Mx^{|\cL^o_i|_M}+1 \eqno(A.1)$$
where the summation is taken over all $M=\{1,\dots,\ell_i\}-I$, where $I$ is a proper sub-interval of $\{1,\dots,\ell_i\}$.  That is, $I$ is nonempty and not the whole interval.  Let us consider the example of an open orbit list $\cL^o_1=\{7,10,8\}$.  Then $\ell_1=3$, and $|\cL^o_1|=7+10+8$, so $T^o_1(x)=x^{25}$.  The proper sub-intervals of $\{1,\dots,\ell_1\}=\{1,2,3\}$ are $I=\{1\},\{2\},\{3\},\{1,2\},\{2,3\}$.  Thus the possibilities for $M=\{1,2,3\}-I$ are $\{2,3\}$, $\{1,3\}$, $\{1,2\}$, $\{3\}$, $\{1\}$.  This gives
$$S^o_1(x)=1+x^{10+8}+x^{7+8}+x^{7+10}+x^8+x^7.$$

If $\cL_i=\cL^c_i$ is an closed orbit list, then we consider $\{1,\dots,\ell_i\}$ to be an interval with a cyclic ordering.  In this case we define 
$$T^c_i(x)=x^{|\cL^c_i|}-1,{\rm\ \ and \ \ }S^c_i(x) = \sum_Mx^{|\cL^c_i|_M}+\ell_i^c \eqno(A.2)$$
where the summation is taken over all $M=\{1,\dots,\ell_i\}-I$, where $I$ is a proper cyclic sub-interval of $\{1,\dots,\ell_i\}$.  That is, $I$ is nonempty and not the whole interval.

\proclaim Theorem A.1. The characteristic polynomial for the matrix representation (4.1--3) is:
$$ \eqalign{ \chi(x) =& (x-d)\prod_{i=1}^{\mu_c}T^{c}_i(x)\prod_{j=1}^{\mu_o}T^{o}_j(x) +
(x-1)\sum_{i=1}^{\mu_c}S^{c}_i\{\prod_{k \neq i}T^{c}_k(x)\prod_{j=1}^{\mu_o}T^{o}_j(x)\} \cr
&+(x-1)\sum_{j=1}^{\mu_o}S^{o}_j\{\prod_{k \neq j}T^{o}_k(x)\prod_{i=1}^{\mu_c}T^{c}_i(x)\}. \cr}\eqno{(A.3)}$$

For each orbit $\cO_{i,j}$, there is an exceptional locus $\Sigma_{i,j}$ and its image $\alpha_{i,j}$.  We have the corresponding subset of ordered basis with $N_{i,j}$ elements,  $\cB_{i,j} = \cF^{i,j}_{N_{i,j}}, \cF^{i,j}_{N_{i,j}-1}, \dots, \cF^{i,j}_1$ where $\cF^{i,j}_k $ is the exceptional fiber $\pi^{-1} p$ in $X$ over $p=f^{N_{i,j}-k}(\alpha_{i,j})$. Let us write down the ordered basis as $$H_X,  \cB^c_{1,\ell^c_1}, \dots, \cB^c_{1,1},\dots, \cB^c_{\mu_c,\ell^c_{\mu_c}}, \dots ,\cB^c_{\mu_c,1},   \cB^o_{1,\ell^o_1}, \dots, \cB^o_{1,1},\dots, \cB^o_{\mu_o,\ell^o_{\mu_o}}, \dots ,\cB^o_{\mu_o,1}.$$ To explain the matrix representation and computation for the characteristic polynomial, let us use  $c^c_{i,j}$ for the column corresponding to the last element of $\cB^c_{i,j}$ and $c^o_{i,j}$ for the column corresponding to the last element of $\cB^o_{i,j}$. We also use $r^c_{i,j}$ for the row corresponding to the first element of $\cB^c_{i,j}$ and $r^o_{i,j}$ for the row corresponding to the first element of $\cB^o_{i,j}$. Using the ordered basis the resulting matrix is 
$$ M=\pmatrix{  d         & B^c_1 &\dots  &B^c_{\mu_c}& B^o_1&\dots  & B^o_{\mu^o} \cr
                                  C^c_1& D^c_1& \dots &                       &             & \dots & \cr
                                  \vdots &            & \ddots &                     &             &           & \cr
                                  C^c_{\mu_c}& &            &D^c_{\mu_c}&             & \dots & \cr
                                  C^o_1&            & \dots &                         &D^o_1&            &\cr
                                  \vdots &             &          &                         &            & \ddots &\cr
                                  C^o_{\mu_o}& &           &                        &           &             & D^o_{\mu_o}\cr}$$
where $B^c_i$ is the $ |\cL^c_i|\times 1$ matrix where $1$ marks
the end of each subset $\cB^c_{i,j}$ of the basis, and thus
 $$ B^c_i = [ 0,\dots,0,1;0,\dots,0,1; 0, \dots,0,1; \dots]$$ 
 and $C^c_i$ is the $ 1 \times |\cL^c_i| $ matrix with $1-d$
marking the beginning of each subset $\cB^c_{i,j}$ of the basis
$\cB^c_{i,\ell^c_i}, \dots ,\cB^c_{i,1}$.  Thus
  $$ C^c_i = [1-d,0, \dots, 0; 1-d,0, \dots, 0; \dots]^t.$$ 
  $B^o_i$ and $C^o_i$ are constructed similarly.  $D^c_i$ is the
$|\cL^c_i|\times |\cL^c_i|$ matrix with 1 in the lower
off-diagonal positions and $0$ for $(c^c_{i,j}, r^c_{i,j+1})$ for
$j=1,\dots ,\ell^c_{i}-1$ and $(c^c_{i,\ell^c_i},r^c_{i,1})$
entries. $D^o_i$ is the $|\cL^o_i|\times |\cL^o_i|$ matrix with
$-1$ in $(c^c_{i,\ell^c_i},r^c_{i,1})$ position and 1 in the lower
off-diagonal positions and $0$ in $(c^o_{i,j}, r^o_{i,j+1})$ for
$j=1,\dots,\ell^o_{i}-1$ positions. Otherwise, $(c^c_{i,j},
r^c_{k,t}), (c^c_{i,j}, r^o_{k,t}), (c^o_{i,j}, r^c_{k,t}),$ and
$(c^o_{i,j}, r^o_{k,t})$ have $-1$ for every remaining combination
of $i,j,k,$ and $t$. All other points are zeros. 

The characteristic polynomial is $\chi(x) = \pm {\rm det} (M-xI)$ with the sign chosen to make it monic.  Now we perform some row  operations on $M-xI$. The top row of $M_xI$ is  $[d-x, B^c_1, \dots, B^c_{\mu_c}, \dots ]$. We add this to every row $r^c_{i,j}$ and $r^o_{i,j}$. 
This produces the matrix
$$ \tilde M=\pmatrix{  d-x      & B^c_1 &\dots  &B^c_{\mu_c}& B^o_1&\dots  & B^o_{\mu^o} \cr
                                  \tilde C^c_1& A^c_1&  &                       &             & \ & \cr
                                  \vdots &            & \ddots &                     &             &           & \cr
                                  \tilde C^c_{\mu_c}& &            &A^c_{\mu_c}&             & & \cr
                                  \tilde C^o_1&            &  &                         &A^o_1&            &\cr
                                  \vdots &             &          &                         &            & \ddots &\cr
                                  \tilde C^o_{\mu_o}& &           &                        &           &             & A^o_{\mu_o}\cr} \eqno{(\tilde M)}$$ where  $\tilde C^c_i$ ,$\tilde C^o_i$ are the same as $C^c_i$, $C^o_i$ except that the entry $1-d$ is now changed to $1-x$. The new diagonal blocks $A^c_i$ have $-x$ in the diagonal, 1 in the lower off-diagonal and 1 in the right hand column of the first row, and $A^o_i$ have $-x$ in the diagonal and 1 in the lower off-diagonal.  That is
$$A^c_i = \left[\matrix { -x&&&1\cr
1&-x&&\cr&\ddots&\ddots&\cr &&1&-x   } \right]_{|\cL^c_i| \times |\cL^c_i|},\ \ A^o_i = \left[\matrix { -x&&&\cr
1&-x&&\cr&\ddots&\ddots&\cr &&1&-x   } \right]_{|\cL^o_i|\times |\cL^o_i|}.$$ The main thing we have accomplished is that except for the first row and column, our matrix $(\tilde M)$ has block diagonal form with blocks $A^c_i$'s, and $A^o_i$'s. The determinant of the each block  has the following form: 
\medskip
\proclaim Lemma A.2. $$ T^{c}_i(x)=(-1)^{|\cL^c_i|}{\rm det} A^c_i=x^{|\cL^c_i|}-1$$ and $$T^{o}_i(x)=(-1)^{|\cL^o_i|}{\rm det}A^o_i=x^{|\cL^o_i|}.$$ 

\medskip

We will evaluate the determinant of $(\tilde M)$ by expanding in minors, going down the left hand column. $1,1$-minor is already block diagonal matrix, thus we use Lemma $A.2$ to take the determinant:
$$\epsilon \prod_{i=1}^{\mu_c}T^{c}_i(x)\prod_{j=1}^{\mu_o}T^{o}_j(x)$$ where $\epsilon=(-1)^{\sum^{\mu_c}_{i=1}|\cL^c_i|+\sum^{\mu_o}_{i=1}|\cL^o_i|}$. Now we consider the $(r_{i,j}, 1)$-minor which is obtained by eliminating the first column and the $r_{i,j}$ row. For us it is more convenient to move the first row to the $r_{i,j}$ position by interchanging two rows. The resulting matrix is almost block diagonal except the $r_{i,j}$ row and its determinant is $(-1)$ times the determinant of the corresponding minor. Let us denote $\hat A_i(j)$ for the diagonal block obtained from  $A_i$ by replacing the row $r_{i,j}$ by $B_i$. For a fixed $i$, $r_{i,j},1$-minor for all $1 \le j \le \ell_i$ correspond to $\cL_i$. The sum of determinant of  $\hat A_i(j)$ for $1 \le j \le \ell_i$ is following:
\medskip
                                        
\proclaim Lemma A.3. $$S^{c}_i(x) =(-1)^{|\cL^c_i|+1} \sum^{\ell^c_i}_{j=1} {\rm det} \, \hat A^c_i (j)$$and $$S^{o}_i(x) = (-1)^{|\cL^o_i|+1} \sum^{\ell^o_i}_{j=1} {\rm det} \, \hat A^o_i (j).$$         
                                  
\noindent {\it Proof.} Let us consider a closed orbit list $\cL^c =\{N_1,\dots, N_\ell\}$. We will compute the determinant of $A^c(j)$ by expanding in minors, going to the right the $r_j$ row. Notice that the $1$'s in $r_j$ row, are not on the diagonal. The minor corresponding the $k$-th $1$ on the left hand side of the diagonal, i.e. $k < j$, is $x^{N_k + \dots N_{j-1}}$. The minor corresponding the $k$-th $1$ on the left hand side of the diagonal, i.e. $i \ge j$, is $x^{N_1 + \dots N_{k-1} + N_{j+1}+ \dots N_\ell}$.  Summing all these monomials for all $k$'s, we have the determinant of $A^c_i(j)$. For the open orbit list, all the minor corresponding the $k$-th $1$ on the left hand side of the diagonal, i.e. $k < j$, are zero since the corresponding matrices are lower triangular matrix with zeros on the diagonal.   

\medskip

For the $(r^c_{i,j}, 1)$-minor, $1$'s in the off-diagonal block would produce zero block on the diagonal. For the $1$'s on the diagonal block will keep other diagonal blocks as diagonal blocks. Thus by Lemma $A.2$ and Lemma $A.3$, sum of the determinant of the  $(r^c_{i,j}, 1)$-minor for $j= 1,\dots, \ell_i$ is given by  $$\epsilon  S^c_i(x)\prod_{k\ne i}T^{c}_k(x)\prod_{j=1}^{\mu_o}T^{o}_j(x)$$ where $\epsilon=(-1)^{\sum^{\mu_c}_{i=1}|\cL^c_i|+\sum^{\mu_o}_{i=1}|\cL^o_i|}$.  

\medskip

\noindent{\it Proof of Theorem A.1.} Combining the previous arguments completes the proof.

\bigskip
\centerline{\bf References}
\medskip

\item{[AABHM]}  N. Abarenkova, J.-C. Angl\`es d'Auriac, S. Boukraa, S.
Hassani, and J.-M. Maillard, Rational dynamical zeta functions for
birational transformations, Phys.\ A 264 (1999), 264--293.

\item{[AABM]} N. Abarenkova, J.-C. Angl\`es d'Auriac, S. Boukraa,  and
J.-M. Maillard,  Growth-complexity spectrum of some discrete dynamical
systems, Physica D 130  (1999), 27--42.

\item{[B]} E. Bedford, On the dynamics of birational mappings of the plane,
J. Korean Math. Soc. 40 (2003), 373--390.

\item{[BMV]} M.P. Bellon, J.-M. Maillard, and C.-M. Viallet, Integrable Coxeter groups, Phys. Lett. A 159 (1991), 221--232.

\item{[BV]} M.P. Bellon and C.-M. Viallet, Algebraic entropy, Comm.\ Math.\
Phys.\ 204 (1999), 425--437.

\item{[BTR]} M. Bernardo, T.T. Truong and G. Rollet, The discrete
Painlev\'e I equations: transcendental integrability and asymptotic
solutions,  J. Phys. A: Math. Gen., 34 (2001), 3215--3252.

\item{[BHM]} S. Boukraa, S. Hassani, J.-M. Maillard,  Noetherian mappings,
Physica D, 185 (2003), no. 1, 3--44. 

\item{[BM]}  S. Boukraa and J.-M. Maillard, Factorization properties of
birational mappings, Physica A 220 (1995), 403--470.

\item{[DF]}  J. Diller and C. Favre, Dynamics of bimeromorphic maps of
surfaces, Amer. J. of Math., 123 (2001), 1135--1169.

\item{[DS]} T.-C. Dinh and N. Sibony, Une borne sup\'erieure pour
l'entropie topologique d'une application rationnelle.  arXiv:math.DS/0303271

\item{[FV]} G. Falqui and C.-M. Viallet, Singularity, complexity, and quasi-integrability of rational mappings, Comm. Math. Phys. 154 (1993), 111--125.

\item{[FS]} J-E Forn\ae ss and N. Sibony, Complex dynamics in higher
dimension: II, Annals of Math. Stud., vol. 137, Princeton University Press,
1995, pp. 135--182.

\item{[G]} V. Guedj, Dynamics of polynomial mappings of ${\bf C}^2$,
Amer. J. of Math.,  124 (2002), 75--106.

\item{[HP]} F.R. Harvey and J. Polking, Extending analytic objects, Comm. Pure Appl. Math. 28 (1975), 701--727.

\item{[RGMR]}  A. Ramani, B. Grammaticos, J.-M. Maillard and G. Rollet,
Integrable mappings from matrix transformations and their singularity
properties, J. Phys. A: Math. Gen. 27 (1994), 7597--7613.

\item{[R1]} K.V. Rerikh, Cremona transformation and general solution of one dynamical system of the static model, Physica D, 57 (1992), 337--354.

\item{[R2]} K.V. Rerikh, Non-algebraic integrability of the Chew-Low reversible dynamical system of the Cremona type and the relation with the 7th Hilbert problem (non-resonant case), Physica D, 82 (1995), 60--78.

\item{[R3]} K.V. Rerikh, Nonalgebraic integrability of one reversible dynamical system of the Cremona type, J. of Math. Physics, 39 (1998), 2821--2832.

\item{[RS]} A. Russakovskii and B. Shiffman, Value distribution of
sequences of rational mappings and complex dynamics, Indiana U. Math.\ J.,
46 (1997), 897--932.

\item{[S]} N. Sibony,  Dynamique des applications rationnelles de ${\bf
P}^k$, Pano.\ Synth.\ 8 (1999), 97--185.

\item{[V]} C. Viallet, On some rational Coxeter groups, Centre de Recherches Math\'ematiques, CRM
Proceedings and Lecture Notes, Volume 9, 1996, 377--388.

\bigskip

\rightline{Department of Mathematics}

\rightline{Indiana University}

\rightline{Bloomington, IN 47405}

\rightline{\tt bedford@indiana.edu}
\bigskip
\rightline{Department of Mathematics}

\rightline{Syracuse University}

\rightline{Syracuse, NY 13244}

\rightline{\tt kkim26@syr.edu}
\end